\documentclass[a4paper]{article}
\pdfoutput=1 
\usepackage{graphicx, lmodern, geometry, xcolor}
\usepackage[utf8]{inputenc}
\usepackage[T1]{fontenc}
\usepackage[british]{babel}
\usepackage{mathtools, amsopn, amssymb, amsfonts, amsthm}
\usepackage[algoruled,noline,noend]{algorithm2e}

\usepackage{enumitem, booktabs, float} 

\usepackage{hyperref} 

\setlist{noitemsep} 
\setlist[enumerate,1]{label = \roman*)}
\hypersetup{
    breaklinks = {true}, 
    pdftitle = {Computations with necklaces on elliptic curves},
    pdfauthor = {Marusia Rebolledo and Christian Wuthrich},
}

\DeclareEmphSequence{\bfseries\color{darkblue}}   

\definecolor{darkblue}{rgb}{0,0,0.7}
\newtheoremstyle{mythm}{}{}{\color{darkblue}}{}{\bfseries\color{darkblue}}{.}{ }{}
\theoremstyle{mythm}
\newtheorem{thm}{Theorem}
\newtheorem{lem}[thm]{Lemma}
\newtheorem{prop}[thm]{Proposition}

\newtheoremstyle{mydef}{}{}{}{}{\bfseries}{.}{ }{}
\theoremstyle{mydef}
\newtheorem*{defn}{Definition}

\newtheorem*{example}{Example}
\newtheorem*{rmk}{Remark}

\newcommand{\QQ}{\mathbb{Q}}
\newcommand{\ZZ}{\mathbb{Z}}
\newcommand{\FF}{\mathbb{F}}
\newcommand{\PP}{\mathbb{P}}

\DeclareMathOperator{\GL}{GL}
\DeclareMathOperator{\Gal}{Gal}
\DeclareMathOperator{\Aut}{Aut}
\DeclareMathOperator{\End}{End}
\DeclareMathOperator{\Frob}{Fr}
\DeclareMathOperator{\PGL}{PGL}
\DeclareMathOperator{\Tr}{Tr}
\DeclareMathOperator{\Norm}{N}
\DeclareMathOperator{\red}{red}
\DeclareMathOperator{\Stab}{Stab}

\newcommand{\vv}{\mathfrak{v}}
\newcommand{\cmset}{\mathcal{CM}}
\renewcommand{\geq}{\geqslant}
\renewcommand{\leq}{\leqslant}

\newcommand\nrfootnote[1]{%
  \begin{NoHyper}
  \renewcommand\thefootnote{}\footnote{#1}%
  \addtocounter{footnote}{-1}%
  \end{NoHyper}
}

\begin{document}

\title{Computing with necklaces on elliptic curves}
\author{Marusia Rebolledo and Christian Wuthrich}
\maketitle

\begin{abstract}
We present computational algorithms to work with points on the modular curve associated to the normaliser of a non-split Cartan group of prime level~$p$.
Rather than working with explicit equations, we represent these points using the moduli interpretation of necklaces in the $p$-torsion of elliptic curves.
We use our methods to investigate for which primes~$\ell\neq p$ two rational points with complex multiplication can have equal reduction modulo~$\ell$.
\nrfootnote{%
The first named author is supported by the ANR Projects ANR-20-CE40-0003 Jinvariant and ANR-23-CE40-0006 GAEC.
The second named author is partially supported by the Engineering and Physical Sciences Research Council, Grant UKRI071.
}
\end{abstract}


\section{Introduction}

Concrete calculations with modular curves have attracted a lot of attention in recent years.
Specifically, the study of algorithms working with explicit isogenies and torsion points on elliptic curves both over number fields and finite fields are an active area.
These correspond to points on the modular curves $X_0(N)$ and $X_1(N)$.
The present article grew out of the attempt to work explicitly with points on a different modular curve, namely the one associated to the normaliser of a non-split Cartan subgroup.
In~\cite{necklaces}, the authors have introduced a moduli interpretation for this curve, which is the starting point for the algorithms presented here.

Let $p$ be an odd prime.
We denote by $X$ the modular curve $X_{\mathrm{nsp}}^+(p)$, which is a smooth projective curve defined over the rational numbers.
A model can be obtained by taking the quotient of the modular curve $X(p)$ by the normaliser of a choice of a non-split Cartan subgroup in~$\GL_2(\FF_p)$.
Let $Y=Y_{\mathrm{nsp}}^+(p)$ be the affine curve obtained by omitting the cusps.
For a field of characteristic different from $p$, the points in $Y(\bar k)$ can be viewed as $\bar{k}$-isomorphism classes of pairs $(E,\vv)$ where $E$ is an elliptic curve defined over $\bar k$ and $\vv$ is a \emph{necklace} in~$E[p]$.
Here a necklace is a particular arrangement of the $p+1$ distinct cyclic subgroup of order $p$ in~$E$.
The precise definition is given in~\cite{necklaces} and repeated in Section~\ref{necklace_background_subsec}.

Equations giving (possibly singular) models for this curve are currently only known for prime level up to $23$ by the works of several authors~\cite{baran, halberstadt, ligozat, mercuri_schoof}.
They are listed together with plenty of other information on the \textsc{lmfdb} data base~\cite{lmfdb}.
Although this is not the focus of the article, the study of the curve
$X=X_{\mathrm{nsp}}^+(p)$ is motivated by Serre’s uniformity conjecture~\cite{serre72}, which could be resolved if the $\QQ$-rational points of $X$ are shown to consist solely of \textsc{cm}~points.
The approach using explicit equations has allowed to determine $X(\QQ)$ in the cases of level $13$ and $17$ using the quadratic Chabauty method in~\cite{cursed, cursed2}.

In this article, we will avoid the use of equations for $X$ and instead,
we make possible concrete calculations using our moduli description of $X$ introduced in \cite{necklaces}.
One way to give a concrete description of a necklace is by listing polynomials $f_0$, $f_1$, \dots, $f_p$ defining the cyclic subgroup schemes $C_0$, $C_1$, \dots, $C_p$ in the order they appear in the necklace.
Even if the necklace is defined over $k$, these polynomials will have coefficients in an extension $L$ which we call the $p$-isogeny field, the smallest extension over which all isogenies of degree $p$ leaving $E$ are defined.
If their codomains, the elliptic curves $E/C_k$, have distinct $j$-invariants, we can also just give these $j$-invariants $j_0$, $j_1$, \dots, $j_p$ as elements in~$L$.
If $p$ is relatively small and $k=\QQ$ or if $k$ is a finite field, we can calculate $L$ and determine all $f_k$; however, for larger $p$ this turns out to be far from efficient.
The use of the $p$-isogeny field compared to the larger field $k\bigl(E[p]\bigr)$ is a good gain: In our case the extension $L/k$ has typically degree $2(p+1)$, which is much smaller than~$2(p^2-1)$.

We present here a first algorithm for $k$ being a number field, which works for any elliptic curve $E$ with a unique $k$-rational necklace~$\vv$.
The implementation in Sage~\cite{sagemath} for $k=\QQ$ can be found at~\cite{implementation}.
The ordering of the subgroups $C_i$ is found using a Frobenius element at a suitable auxiliary prime ideal $\mathfrak{Q}$ in~$L$.
For any prime ideal $\ell$ in $k$, we can then calculate the reduction of $(E,\vv)$ modulo $\ell$ to obtain a representation $(\tilde E,\tilde\vv)$ of a point in $X(\FF_{\ell})$ where $\FF_{\ell}$ is the residue field at~$\ell$.
We have to emphasise that this is really a global problem; we cannot work in a reduction or in a completion as there the curve will usually have more than one necklace and we cannot guess which one is the reduction of a $k$-rational necklace we are after.
Unfortunately, this first algorithm is only really practical for small primes~$p$ and for them we know good models for~$X$.

It is expected that for $p>11$, the only points in $Y(\QQ)$ are given by $(E,\vv)$ with $E$ having complex multiplication.
Therefore, we present a second faster algorithm for such points called \emph{CM points}.
The algorithm discussed in Section~\ref{cm_alg_subsec} calculates the reduced point $(\tilde E,\tilde \vv)\in X(\FF_{\ell})$ directly without having to determine $L$ globally.
It only needs to work out the $p$-isogeny field for the curve $\tilde E$ over the finite field~$\FF_{\ell}$.
We do not know a similar construction for curves without complex multiplication.

As an application of this algorithm, we can make predictions about the following question.
For which distinct primes $p,\ell\geq 5$ do there exist two \textsc{cm}~points $x_1$ and $x_2$ in $X(\QQ)$ such that their reductions in $X(\FF_{\ell})$ are equal?
For instance, when using an explicit Chabauty methods as in~\cite{cursed, cursed2} to determine the rational point on $X$ for a particular~$p$, the auxiliary prime $\ell$ is best chosen such that all rational points have distinct reduction, since this results in $\ell$-adic analytic functions with a single zero in a residue disc.
One could hope that the moduli interpretation will help the use of the geometric Chabauty method as in~\cite{geometric_chabauty} for larger~$p$ without having to rely on explicit equations.
This is one of the main motivations of this paper.

When $p$ increases, the number of points in $X(\FF_{\ell})$ increases quite quickly, which means that we do not expect equal reduction among the finitely many \textsc{cm}~points $x_1$ and $x_2$ in $X(\QQ)$ to happen when $p$ is large.
Having tested all $5 \leq p < 50$ and all $\ell>3$, we found equal reductions only for $p=5$ and $p=7$ and for~$\ell\leq 17$.
The eight such pairs $(x_1,x_2)$ we found are listed in Proposition~\ref{inj_prop}; it is likely that these are the only examples that exist.

This question can be seen as the first step in a deeper investigation into the arithmetic intersection number of points on the curve $X$ analogous to the work of Gross--Zagier in~\cite{gross_zagier} for the modular curve $X_0(p)$.
Recently, this has been studied by Love, Studnia, and Vonk in~\cite{love_studnia_vonk}.
In particular, their work gives a theoretical understanding of these intersections and verifies indeed that our eight pairs are the only such cases of equal reduction among \textsc{cm}~points.

The algorithms in this paper are to our knowledge the first methods found to do explicit calculations with points on $X_{\mathrm{nsp}}^+(p)$ without the use of equations defining the curve.
They have their limitations, but do allow for experimentations on these curves.
Of course, it would be very interesting to find faster or more general methods.
At this stage, we do not know if there are any applications of these algorithms over finite fields, like the theory of isogeny volcanoes~\cite{sutherland_volcanoes} has.

The paper is structured as follows.
Section~\ref{background_sec} recalls the definition of necklaces and reviews background results used later.
The algorithm to calculate any necklace over a number field is explained in Section~\ref{algo1_sec}.
The reduction of necklaces, both for general as well as \textsc{cm}~points is contained in Section~\ref{red_sec}, while Section~\ref{comparing_sec} considers the question of when two points have the same reduction.
The short appendix~\ref{nu_pts_app} gives a table that allows to calculate $\# X(\FF_{\ell})$.

\paragraph{Acknowledgements.} We are very grateful to the referee for carefully reading our article and for making many suggestions that have improved the paper.

\section{Background}\label{background_sec}

In this section we recall the moduli interpretation of the modular curve $X=X_{\mathrm{nsp}}^+(p)$ described in~\cite{necklaces} and state some preliminary results.
We wish to point to~\cite{serre} for basic facts about the Galois representation $E[p]$ and the notions of Cartan subgroups and their normalisers.

\subsection{Necklaces}\label{necklace_background_subsec}
Let $p$ be an odd prime.
We fix throughout the article a generator~$\gamma$ of the cyclic group~$\FF_{p^2}^{\times}$.
By a \emph{necklace} on an elliptic curve, we will understand a non-oriented $\gamma$-necklace as defined in~\cite{necklaces} whose definition we are going to recall now.

\begin{defn}
  Let $E$ be an elliptic curve defined over a field $k$ whose characteristic is different from~$p$.
  A \emph{necklace} in the $p$-torsion of $E$ is defined to be an equivalence class $\vv$ of an ordering $(C_0,C_1,\dots,C_p)$ of all cyclic subgroups of order $p$ in $E(\bar{k})$ satisfying the condition that there is an element $h\in \PGL\bigl(E[p]\bigr)$ with $h(C_i) = C_{i+1 \bmod{p+1}}$ and such that there is a matrix in $h$ whose characteristic polynomial is equal to the minimal polynomial of~$\gamma$.
  Two such lists are equivalent if one can be obtained from the other by a cyclic permutation and a reversal $w\colon (C_0, C_1, \dots, C_p) \mapsto (C_p,\dots, C_1,C_0)$ if needed.
\end{defn}

We view the cyclic subgroups $C_i$ as pearls and the necklace as placing these pearls on a regular $(p+1)$-gon forming the picture of a pearl necklace:

\begin{center}
  \includegraphics[height=3cm]{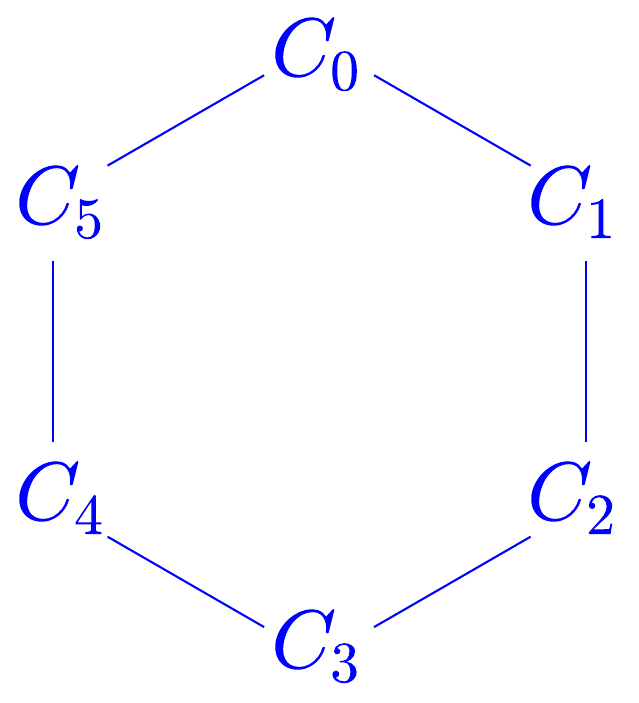}
\end{center}

Turning a necklace $\vv$ (acting by $h$) or flipping it from one side to the other (acting by the involution $w$) does not change the necklace by the equivalence introduced above.
The subgroup of $\PGL(E[p])$ stabilising $\vv$ is generated by $h$ and $w$: this is the normaliser of the non-split Cartan subgroup generated by $h$.
Conversely, if $N$ is the normaliser of a non-split Cartan subgroup $C$ and $h$ an element that generates $C$ with characteristic polynomial equal to the minimal polynomial of $\gamma$, then acting successively by $h$ on a given cyclic subgroup $C_0$ of order $p$ gives the only necklace stabilised by~$N$.
There are $p(p-1)/2$ necklaces (see \cite[Corollary 3]{necklaces}).

As the terminology is a little cumbersome to repeat often, we will introduce the neologism \emph{nonoca} as an abbreviation for ``normaliser of a non-split Cartan subgroup'' in $\PGL\bigl(E[p]\bigr)$.
Recall that any nonoca is isomorphic to a dihedral group of order $2(p+1)$.

There is a characterisation of necklaces without reference to $h$:
Let $\xi = t^2/(t^2-n)$ where $t=\Tr(\gamma)$ and $n=\Norm(\gamma)$.
A list $(C_0, C_1,\dots, C_p)$ represents a necklace if the cross-ratio $[C_i,C_{i+1};C_{i+2},C_{i+3}]$ of any four consecutive pearls is equal to $\xi$.
See Proposition~5 in~\cite{necklaces}.

Any choice of three pearls $C_0$, $C_1$ and $C_2$ yields a unique necklace such that $C_1$ is adjacent to $C_0$ and~$C_2$.
See Lemma~2 in~\cite{necklaces}.

The absolute Galois group $G_k$ of $k$ acts on the cyclic subgroups of order $p$ of $E$; we write $\rho\colon G_k \to \PGL\bigl(E[p]\bigr)$ for the corresponding map.
This induces an action of $G_k$ on the set of necklaces on~$E$.

Let $X = X_{\text{nsp}}^{+}(p)$ be the modular curve associated to the normaliser of a non-split Cartan subgroup of~$\GL_2(\FF_p)$.
The affine curve obtained by removing the cusps is denoted by $Y=Y_{\text{nsp}}^{+}(p)$.
The following is proved in Section~2.3 in~\cite{necklaces}.

\begin{prop}
  Suppose that the characteristic of $k$ is different from $p$.
  There is a bijection between points on $Y(\bar{k})$ and $\bar{k}$-isomorphism classes of pairs $(E,\vv)$ where $E$ is an elliptic curve defined over $\bar{k}$ and $\vv$ is a necklace in the $p$-torsion of $E$.
  Any point in $Y(k)$ is represented by $(E,\vv)$ with $E$ defined over $k$ and, in case $j(E)\not\in \{0,1728\}$, the necklace $\vv$ is also defined over $k$.
\end{prop}

For a point in $Y$ with $j$-invariant different from $0$ and $1728$ to be $k$-rational is equivalent to asking that $\rho(G_k)$ is contained in a nonoca, namely the nonoca stabilising the corresponding necklace.
The situation is more subtle for the points with $j$-invariant in $\{0,1728\}$ as explained in Section~\ref{extra_aut_subsec}.

Another important notion that we will be using is the following.
Two pearls $C_i$ and $C_j$ are \emph{antipodal} in a necklace $\vv=(C_0,C_1,\dots,C_p)$ if $i\equiv j+(p+1)/2 \pmod{p+1}$.
That is they are diametrically opposed when the necklace is represented as a regular $(p+1)$-gon:

\begin{center}
  \includegraphics[height=3cm]{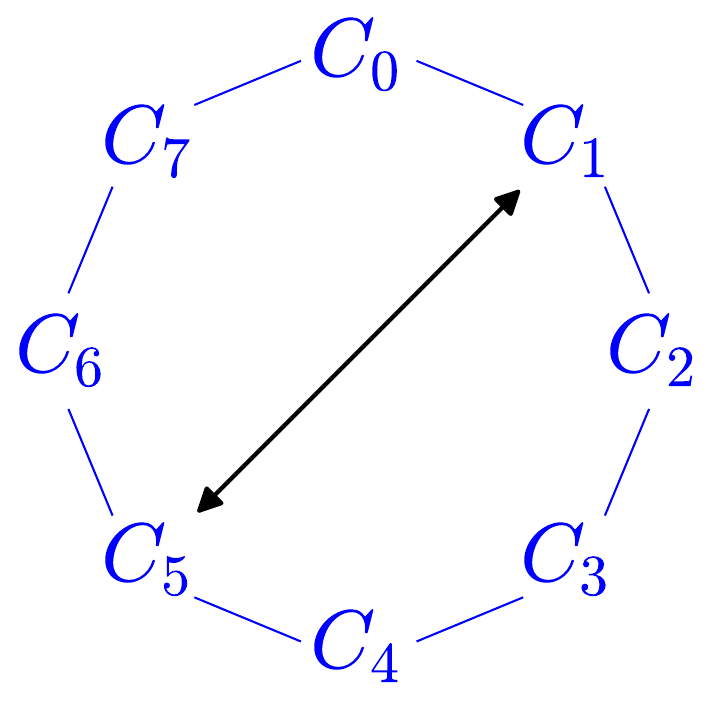}
\end{center}

\subsection{The number of rational necklaces on a given elliptic curve}\label{number_subsec}
We aim to describe the number of $k$-rational necklaces on an elliptic curve defined over a field~$k$.
For instance, we will see that for $k=\QQ$ there is usually at most one such necklace.

\begin{lem}\label{nb_gen_lem}
  Let $E$ be an elliptic curve over a field $k$ whose characteristic is coprime to~$p$.
  The number of necklaces on $E$ defined over $k$ could be either $0$, $1$, $2$ (only if $p\equiv 1\pmod{4}$), $3$ (only if $p\equiv 3 \pmod{4}$), $(p-1)/2$, $(p+3)/2$ or all $p(p-1)/2$ of them.
\end{lem}

Denote by $G\leq \PGL\bigl(E[p]\bigr)$ the image of the Galois representation $\rho\colon G_k \to \PGL\bigl(E[p]\bigr)$.
Given an element $g\in \PGL\bigl(E[p]\bigr)$ represented by an element $M\in \GL\bigl(E[p]\bigr)$, we will write
\[
  \delta(g) = \Bigl(\frac{\Tr(M)^2 - 4 \det(M)}{p}\Bigr) \in \{-1,0,1\},
\]
where $\bigl(\tfrac{\cdot}{\cdot}\bigr)$ denotes the Legendre symbol.
If $\delta(g)\neq 0$, then $g$ lies in a unique Cartan subgroup.
The element $g$ belongs to a split Cartan subgroup precisely when $\delta(g)=1$, in which case we will say that $g$ is \emph{split}.
Otherwise, when $\delta(g)=-1$, we will say that $g$ is \emph{non-split}.
Recall that a non-trivial element of $\PGL(E[p])$ is split if and only if it fixes two subgroups of $E[p]$.

\begin{proof}
  A necklace $\vv$ is defined over $k$ if and only if $G$ is contained in the nonoca stabilising~$\vv$.
  Hence there is no necklace defined over $k$ on $E[p]$ in the case when $G$ is not contained in any nonoca and there is a unique one in the case when $G$ is contained in a single nonoca.
  If $G$ is trivial, then all $p(p-1)/2$ necklaces are defined over $k$.
  We may assume now that $G$ is non-trivial.

  Suppose now that $E$ has at least two necklaces $\vv$ and $\vv'$ defined over $k$ and hence $G\subset N\cap N'$ where $N$ and $N'$ are the nonocas stabilising $\vv$ and $\vv'$, respectively.
  As any two non-split Cartan subgroups intersect trivially, the only possible non-trivial elements in the intersection of two nonocas $N$ and $N'$ have order~$2$.

  In what follows, we will say that an element of $G$ \emph{fixes} or \emph{flips} a necklace $\vv$ if it fixes or flips the associated oriented necklace as defined in~\cite{necklaces}.

  If a non-trivial element $g\in G$ is split, it is not in any non-split Cartan subgroup, hence any necklace defined over $k$ is flipped by $g$, that is to say $g$ appears as a reflection.
  It follows that the subgroups $A$ and $B$ fixed by $g$ are antipodal in the necklace and $g$ is the reflection whose axis is the line through $A$ and~$B$.
  By Lemma~8 in~\cite{necklaces} there are $(p-1)/2$ such necklaces.

  If $g$ is non-split, it is in a unique non-split Cartan subgroup~$C_g$ and hence it fixes a unique necklace, namely the necklace stabilised by the normalizer of~$C_g$.
  In this necklace, $g$ appears as the rotation of angle~$\pi$.

  We count the number of necklaces that $g$ flips as follows.
  The non-split elements of order $2$ form a single conjugacy class in $\PGL_2(E[p])$ of cardinality $p(p-1)/2$; every non-split Cartan subgroup contains exactly one such non-split element of order~$2$.
  Therefore, each of them flips the same number of necklaces.
  There are $(p+1)/2$ reflections of a regular $(p+1)$-gon that fix no corner.
  Therefore $g$, and any of its conjugate elements, will flip $(p+1)/2$ necklaces.
  In total, the non-split $g$ belongs to $1+(p+1)/2=(p+3)/2$ nonocas.

  This concludes the case when $G$ is cyclic generated by $g$: We obtain $(p-1)/2$ necklaces defined over $k$ if $g$ is split, or $(p+3)/2$ necklaces over $k$ if $g$ is non-split.

  Suppose now that $G\subset N\cap N'$ contains two distinct non-trivial elements $g_1$ and~$g_2$.
  In fact, $G$ is then isomorphic to the abelian, but non-cyclic group of order~$4$ because these are the largest subgroups of a nonoca containing no elements of larger order.
  We claim that of the three elements $g_1$, $g_2$, and $g_3=g_1g_2$ at most one is split.
  Indeed, if there were two split elements, the two distinct necklaces $\vv$ and $\vv'$ would have two distinct antipodal pairs of pearls in common, which is impossible by Lemma~8 in~\cite{necklaces}.

  Hence we may suppose that $g_1$ and $g_2$ are non-split.
  As elements of order~$2$ come from elements in $\GL\bigl(E[p]\bigr)$ with trace zero, we find $\delta(g_3) = \bigl(\tfrac{-1}{p}\bigr)\,\delta(g_1)\,\delta(g_2) = \bigl(\tfrac{-1}{p}\bigr)$.
  Thus $g_3$ is split if $p\equiv 1\pmod 4$ and non-split if $p\equiv 3\pmod 4$.

  We claim that any necklace $\vv$ fixed by a non-split element $g\in G$ is automatically flipped by any other element $g'\in G$ using that they commute.
  Indeed, for any $A\in \PP(E[p])$, $A$ and $g(A)$ are antipodal in $\vv$, and so are $g'(A)$ and $gg'(A)$.
  In the necklace $g'(\vv)$, it follows that $g'(A)$ and $g'g(A)=gg'(A)$ are antipodal and so are $A=g'g'(A)$ and $g(A)=g'g'g(A)$.
  However, by Lemma~8 in~\cite{necklaces} there is at most one necklace with a fixed pair of antipodals.
  Hence $g'(\vv)=\vv$ and $g'$ must flip $\vv$ as $g$ already fixes it.

  First the case $p\equiv 1\pmod 4$:
  A necklace $\vv$ defined over $k$ is flipped by $g_3$ since $g_3$ is split.
  If it were also flipped by $g_1$, then it is fixed by~$g_2=g_1g_3$.
  Thus exactly one of $g_1$ or $g_2$ must fix~$\vv$.
  We deduce that there are exactly two necklaces defined over $k$: one fixed by $g_1$, which is flipped by $g_3$ and $g_2$, and one fixed by $g_2$ and flipped by $g_1$ and~$g_3$.

  Finally the case $p\equiv 3\pmod 4$:
  With the same argument as above, a necklace defined over~$k$ must be fixed by exactly one of $g_1,g_2,g_3$ and then, by the above claim, it is automatically flipped by the others.
  We conclude that there are exactly three necklaces defined over~$k$.
\end{proof}

\begin{lem}\label{number_of_necklaces_lem}
Let $E$ be an elliptic curve defined over a finite field $\FF$ of cardinal $\ell^r$ with $\ell$ a prime distinct from~$p$ and $r\geq 1$ an integer.
Let $a$ be the trace of Frobenius such that $\# E(\FF) = \ell^r+1-a$ and set $\delta=\bigl(\tfrac{a^2-4\ell^r}{p}\bigr)\in\{-1,0,1\}$.
Denote by $n_{\text{isog}}$ the number of isogenies of degree $p$ defined over $\FF$ leaving from~$E$.
The number $n_{\mathrm{necklaces}}$ of necklaces on $E$ defined over $\FF$ can be read off the following table:
  \begin{center}\setlength{\tabcolsep}{12pt}\begin{tabular}{cccc}\toprule
    $a$ & $\delta$ & $n_{\text{isog}}$ & $n_{\text{necklaces}}$ \\
    \midrule
    $0$ & $+1$ && $(p-1)/2$ \\
    $0$ & $-1$ && $(p+3)/2$ \\
    $\neq 0$ & $+1$ && $0$ \\
    $\neq 0$ & $-1$ && $1$ \\
    $\neq 0$ & $0$ & $p+1$ & $p(p-1)/2$\\
    $\neq 0$ & $0$ & $1$ & $0$ \\
    \bottomrule
  \end{tabular}\end{center}
\end{lem}

This lemma counts the number of $\FF$-rational necklaces for $\FF$ a finite field of characteristic different from~$p$.
For curves with $j$-invariant different from $0$ and $1728$, this is the same as the number of $\FF$-rational points on $X$ with that $j$-invariant.
For the special two $j$-invariants this is more complicated and discussed in Section~\ref{extra_aut_subsec}.
In the case $\FF=\FF_\ell$, the complete table for all $j$-invariants is given in Appendix~\ref{nu_pts_app}, which contains the above as its final six lines.

\begin{proof}
  As in the proof above, let $G$ be the image of the absolute Galois group in $\PGL\bigl(E[p]\bigr)$.
  It is generated by the image $g$ of the Frobenius, whose image in $\GL\bigl(E[p]\bigr)$ is a matrix $M$ with characteristic polynomial $X^2-a\,X+\ell^r \in \FF_p[X]$.
  By definition $\delta(g) = \delta$.
  In the top row of the table, $g$ is a split element of order 2.
  In the second row, it is a non-split element of order 2.
  If we are in the third row, then $g$ is not in any nonoca since it is split  of order $>2$.
  In the fourth row, $g$ is non-split of order larger than $2$; it belongs to a unique nonoca.
  The last two rows are matrices with repeated eigenvalues, they can either be diagonalisable matrices, that is $g=1$ and all necklaces are rational over $\FF$, or non-diagonalisable, in which case they are in no nonoca.
  In the first case, all $p+1$ isogenies $E\to E'$ of degree $p$ are defined over $\FF$, in the second case, there is a unique such isogeny.
\end{proof}

\begin{prop}\label{nboverQ_prop}
  Let $E$ be an elliptic curve defined over $\QQ$ with~$j(E)\notin \{0, 1728\}$.
  If $p>5$, then there is at most one necklace defined over $\QQ$ on~$E$.
  If $p=5$, then there are at most two.
\end{prop}

\begin{proof}
  Suppose first that $E$ does not have complex multiplication.
  Theorem~1.5 in \cite{furio_lombardo} by Furio and Lombardo, improving results of Zywina~\cite{zywina2015} and Le Fourn and Lemos \cite{lefourn_lemos}, shows that for $p>37$, either $\rho\colon G_{\mathbb{Q}}\to \PGL\bigl(E[p]\bigr)$ is surjective or it has image equal to a whole nonoca.
  They also proved in their Theorem 1.6 that for primes $5<p\leq 37$ the image of $\rho$ cannot be a proper subgroup of a nonoca.
  We deduce from this that if $E$ does not have complex multiplication and $p>5$, then there is at most one necklace defined over $\QQ$ on~$E$.

  Suppose now that $E$ has complex multiplication by an order $\mathcal O$ of an imaginary quadratic field $F$, by which we mean $\End_{\bar \QQ}(E)\cong \mathcal O$.
  If the image of $\rho$ is contained in a nonoca, then $p$ is inert in~$\mathcal O$ and $\rho(G_\QQ)$ is the whole nonoca as showed in Proposition~1.14 in~\cite{zywina2015} or in our Lemma~\ref{cm_lem} below.
  Again, we deduce that there is a unique necklace defined over~$\QQ$.

  The image of complex conjugation under $\rho$ is a class of matrices with trace $0$ and determinant $-1$; therefore it is a split element of order $2$.
  For $p=5$, this fact together with the proof of Lemma~\ref{nb_gen_lem} shows that there are at most $\frac{p-1}{2}=2$ necklaces defined over~$\QQ$.
\end{proof}

\begin{example}
  For $p=5$, the image of $\rho$ can be equal to $C_2\times C_2$, a case denoted by $G(p)=G_3$ in Theorem 1.4 in~\cite{zywina2015}.
  In this case, there are exactly two necklaces defined over $\QQ$ on $E$, as shown by the proof of Lemma~\ref{nb_gen_lem}.
  For example, the curve with Cremona label \href{https://beta.lmfdb.org/EllipticCurve/Q/6975/h/1}{6975d1} has two necklaces for $p=5$ defined over~$\QQ$.
  There are infinitely many such examples as the \href{https://beta.lmfdb.org/ModularCurve/Q/5.30.0.b.1/}{modular curve} in question is of genus~$0$.
\end{example}

\begin{example}
  Note that the image of $\rho$ for the curve \href{https://beta.lmfdb.org/EllipticCurve/Q/98/a/3}{98a3} for $p=3$ consists of only one non-split element of order $2$.
  All three necklaces are defined over $\QQ$ on that curve.
  The corresponding \href{https://beta.lmfdb.org/ModularCurve/Q/3.24.0-3.a.1.1/}{modular curve} has genus $0$, so there are infinitely many examples over~$\QQ$.
\end{example}

\subsection{Elliptic curves with complex multiplication}\label{cm_subsec}
For the convenience of the reader, we will recall some background results for elliptic curves with extra endomorphisms.
We will say $E$ has complex multiplication (\textsc{cm}) if the geometric endomorphism ring is not $\ZZ$, even if the endomorphisms are not defined over the base field of the curve.

\begin{lem}\label{cm_classic_lem}
  Let $E$ be an elliptic curve defined over a number field $K$ with complex multiplication by an order $\mathcal O$ with conductor $f$ in an imaginary quadratic field~$F$.
  \begin{enumerate}
    \item\label{case1} If $p$ is split in $\mathcal O$, then the image of $\rho\colon G_K \to \PGL\bigl(E[p]\bigr)$ is contained in the normaliser of a split Cartan subgroup.
    \item\label{case2} If $p$ is inert in $\mathcal O$, then the image of $\rho$ is contained in  a nonoca.
    \item If $p$ ramifies in $\mathcal O$, then the image of $\rho$ is contained in a Borel subgroup, but not in any normaliser of a Cartan subgroup.
  \end{enumerate}
  In case~\ref{case1} and~\ref{case2}, if moreover $p$ does not divide $f$ nor the absolute discriminant $\Delta_{FK}$ of $FK$ and $E$ has good reduction above $p$, then $\rho$ has image the whole Cartan subgroup if $F\subset K$ or its normalizer if $[FK:K]=2$.
\end{lem}

If $p$ is inert in $\mathcal O$, the necklace $\vv^*$ fixed by the non-split Cartan of Lemma~\ref{cm_classic_lem} is defined over $K$ and the point $[(E,\vv^*)]\in X(K)$ is a Heegner point.
(See also Section~\ref{cm_alg_subsec}).

\begin{proof}
  Over $FK$, the action by $G_{FK}$ and the action by the endomorphism ring $\mathcal{O}$ on $E[p]$ commute.
  This shows that the restriction of $\rho$ to $G_{FK}$ maps into the subgroup $C_{\mathcal O}$ in $\PGL\bigl(E[p]\bigr)$ which is the image of $\Aut_{\mathcal{O}}\bigl(E[p]\bigr)\cong (\mathcal{O}/p\mathcal{O})^{\times}$.
  If $p$ splits in $\mathcal{O}$ this is a split Cartan subgroup and, if $p$ is inert, it is a non-split Cartan subgroup.
  If $p=\mathfrak{p}^2$ for an ideal $\mathfrak{p}$ in $\mathcal{O}$, then the Galois action must fix the subgroup $E[\mathfrak{p}]$ inside $E[p]$ which shows that $\rho(G_{FK})$ lies inside a Borel subgroup, however it can not belong to a split Cartan subgroup as there is no other sub-$\mathcal{O}$-module in $E[p]$ of order~$p$.
  Then any element of a Borel subgroup which is not split has order $p$, and hence it is not in a non-split Cartan.

  As $[FK:K]\leq 2$, the subgroup $\rho(G_{FK})$ has index at most $2$ in $\rho(G_K)$, which implies that it is normal.
  Hence that the image of $\rho$ is in the normaliser of~$\rho(G_{FK})$.

  If, in addition, $p$ does not divide $f\Delta_{FK}$ and $E$ has good reduction above $p$, then the natural injection $G_{FK}\hookrightarrow C_{\mathcal O}$ is an isomorphism, as proved in \cite[Proposition~3.3]{campagna-pengo}.
  This gives the second part of the lemma if $F\subset K$.
  If $[FK:K]=2$, it remains to show that $\rho(G_K)$ is the whole normaliser.
  Suppose it is not, hence $\rho(G_K)=\rho(G_{FK})$ and the action of $G_K$ commutes with the action of $\mathcal O$ on $E[p]$.
  This contradicts the fact that the action of a non-trivial element of $\Gal(FK/K)$ and an element of degree~$2$ in $\mathcal O$ anti-commute.
\end{proof}

When $E$ is defined over $\QQ,$ we can say a little more.
The order $\mathcal O$ is one of the thirteen imaginary quadratic orders of class number one listed in the Table~\ref{cm_overQ_tab} below.
We denote by $D=\Delta_F\cdot f^2$ the discriminant of $\mathcal O$.
The elliptic curve $E_{D}$ is one of minimal conductor among the elliptic curves with complex multiplication by $\mathcal O$ and they are listed with their Cremona label here.

\begin{table}[H]
  \centering
  \caption{\textsc{cm} elliptic curves over $\QQ$}\label{cm_overQ_tab}
  \begin{tabular}{lrr}
    \toprule
    $D = \Delta_F\cdot f^2$   &  $j$  &  $E_D$ \\
    \midrule
    $-3$ & $0$ & \href{https://www.lmfdb.org/EllipticCurve/Q/27a3/}{27a3} \\
    $-4$ & $1728$ & \href{https://www.lmfdb.org/EllipticCurve/Q/32a2/}{32a2} \\
    $-7$ & $-3375$ & \href{https://www.lmfdb.org/EllipticCurve/Q/49a1/}{49a1} \\
    $-8$ & $8000$ & \href{https://www.lmfdb.org/EllipticCurve/Q/256a1/}{256a1} \\
    $-11$ & $-32768$ & \href{https://www.lmfdb.org/EllipticCurve/Q/121b1/}{121b1} \\
    $-12=-3\cdot 2^2$ & $54000$ & \href{https://www.lmfdb.org/EllipticCurve/Q/36a2/}{36a2} \\
    $-16=-4\cdot 2^2$ & $287496$ & \href{https://www.lmfdb.org/EllipticCurve/Q/32a3/}{32a3} \\
    $-19$ & $-884736$ & \href{https://www.lmfdb.org/EllipticCurve/Q/361a1/}{361a1} \\
    $-27=-3\cdot 3^2$ & $-12288000$ & \href{https://www.lmfdb.org/EllipticCurve/Q/27a2}{27a2} \\
    $-28=-7\cdot 2^2$ & $16581375$ & \href{https://www.lmfdb.org/EllipticCurve/Q/49a2/}{49a2} \\
    $-43$ & $-884736000$ & \href{https://www.lmfdb.org/EllipticCurve/Q/1849a1/}{1849a1} \\
    $-67$ & $-147197952000$ & \href{https://www.lmfdb.org/EllipticCurve/Q/4489a1/}{4489a1} \\
    $-163$ & $-262537412640768000$ & \href{https://www.lmfdb.org/EllipticCurve/Q/26569a1/}{26569a1} \\
    \bottomrule
  \end{tabular}
\end{table}

\begin{lem}\label{cm_lem}
  Let $E$ be an elliptic curve over $\QQ$ with complex multiplication by an order $\mathcal O$ of an imaginary quadratic field $F$ of discriminant $\Delta_F$ and let $p>2$ be a prime number.
  Assume that $j(E)\notin \{0,1728\}$.
  \begin{enumerate}
    \item If $\bigl(\tfrac{\Delta_F}{p}\bigr)=1$ then the image of $\rho$ is equal to the normaliser of a split Cartan subgroup.
    \item If $\bigl(\tfrac{\Delta_F}{p}\bigr)=-1$ then the image of $\rho$ is equal to a nonoca.
    \item If $p\mid \Delta_F$ then the image of $\rho$ is contained in a Borel, but not in any normaliser of a Cartan subgroup.
  \end{enumerate}
  Hence, if $p>3$, there is a $\QQ$-rational necklace $\vv^*$ on $E$ only in the case $\bigl(\tfrac{\Delta_F}{p}\bigr)=-1$ and it is then unique.
\end{lem}

\begin{proof}
  From Table~\ref{cm_overQ_tab}, we see that $p\nmid \Delta_F$ implies $p\nmid f$ for all odd $p$.
  It follows that the condition on $\bigl(\tfrac{\Delta_F}{p}\bigr)$ corresponds to the splitting behaviour of $p$ in $\mathcal{O}$, so $\rho(G_\QQ)$ is contained in the normaliser of a split Cartan, the normaliser of a non-split Cartan or a Borel as listed in Lemma~\ref{cm_classic_lem}.

  The case $p\mid \Delta_F$ follows from Lemma~\ref{cm_classic_lem}.
  Now, we wish to show that if $p\nmid \Delta_F$, $\rho(G_\QQ)$ is the whole normaliser of the  Cartan subgroup $C_\mathcal O$ (with the notation of the proof of Lemma~\ref{cm_classic_lem}).

  As stated by Stevenhagen (see Theorem~1.4 in~\cite{bourdon_clark}), class field theory implies that the $p$-ray class field of $\mathcal O$  is included in $F(E[p])$ and has Galois group over $F$ isomorphic to $(\mathcal O/p\mathcal O)^\times/ [\mathcal{O}^{\times}]$ where $[\mathcal{O}^\times]$ is the image of $\mathcal{O}^{\times}$ through $\mathcal O\longrightarrow \mathcal O/p\mathcal O$.
  Hence, the index of $\rho(G_F)$ in $C_{\mathcal{O}}$ divides $\lvert\mathcal{O}^\times/\{\pm 1\}\rvert$, which gives $\rho(G_F)=C_{\mathcal{O}}$ since $j(E)\notin\{0,1728\}$.
  We conclude that $\rho(G_\QQ)$ is the whole normaliser of $C_{\mathcal{O}}$ as in the end of the proof of Lemma~\ref{cm_classic_lem}.

  The following is an alternative proof: Since $j(E)\notin\{0,1728\}$, the curve $E$ is a quadratic twist of $E_D$.
  Therefore the image of $\rho$ is isomorphic to the image of the representation associated to $E_D$.
  From Table~\ref{cm_overQ_tab}, we see that $E_D$ has good reduction for all odd primes $p\nmid \Delta_F$.
  Hence the result follows from the second part of Lemma~\ref{cm_classic_lem} for $K=\QQ$.

  Finally, note that the normaliser of a split Cartan subgroup cannot be contained in a nonoca if $p>3$, justifying the final statement of the Lemma.
\end{proof}

Recall that the cases $j\in\{0,1728\}$, are discussed in Section~\ref{extra_aut_subsec}.

\begin{rmk}
 We state and prove Lemmas \ref{cm_classic_lem} and \ref{cm_lem} in the projective setting, which is sufficient for our purposes and slightly simpler, but they also hold for the representation into $\GL(E[p])$ instead of $\PGL(E[p])$.
 The proof of Lemma \ref{cm_classic_lem} remains unchanged, while the $\GL(E[p])$ version of Lemma \ref{cm_lem} follows, for instance, from Theorem 6.3 in \cite{campagna-pengo}.
\end{rmk}

\subsection{Distinct \texorpdfstring{$j$}{j}-invariants}
Since our goal is to represent necklaces algorithmically, we are particularly interested by the case when the $j$-invariants of the curves $E/C$ for all $C$ cyclic subgroup of order $p$ are distinct.
Indeed, in this case we may represent $\vv=(C_0,C_1,\dots,C_p)$ by an ordered list of the $j$-invariants $(j(E/C_0),j(E/C_1),\dots,j(E/C_p)\bigr)$.
We will make use of this in Section~\ref{algo1_sec}.

\begin{lem}\label{distinct_j_lem}
  Let $E$ be an elliptic curve defined over a field $k$, such that $j(E)\not\in\{0,1728\}$.
  Assume that $k$ is of characteristic $0$ or that $k$ is of characteristic~$\ell\neq p$ and $E$ is ordinary.
  There exists two distinct cyclic subgroups $C$ and $C'$ of order $p$ such that $j(E/C)=j(E/C')$ if and only if $\End(E)$ is an order in an imaginary quadratic field in which $p$ splits.
\end{lem}

\begin{proof}
  Suppose that $C$ and $C'$ are two distinct cyclic subgroups of order $p$ in $E$ such that $E/C$ and $E/C'$ are $\bar k$-isomorphic, say by an isomorphism $\iota$.
  If $\varphi$ and $\varphi'$ are choices of isogenies with $\ker(\varphi)=C$ and $\ker(\varphi')=C'$, then the composition $\alpha=\check\varphi \circ \iota\circ \varphi' \colon E\to E/C'\to E/C\to E$ gives an endomorphism on $E$.

  If the ideal $I =(\alpha)$ is equal to $(p)$ in $\End(E)$, then there is a unit $u\in\End(E)$ such that $\alpha = u\, p$.
  As we have excluded that $j(E)$ is $0$ or $1728$, we must have $u=\pm 1$, but this would imply that $\iota\circ\varphi'=\pm \varphi$, which is impossible if $\varphi$ and $\varphi'$ have distinct kernel.

  It follows that the endomorphism ring of $E$ is larger than $\ZZ$.
  Our hypothesis implies that it is isomorphic to an order $\mathcal{O}$ in an imaginary quadratic field.
  The ideal $I$ is such that $\mathcal{O}/I$ is cyclic of order $p^2$.
  This implies that $I$ is the square of a prime ideal of norm $p$ and hence $p$ splits in $\mathcal{O}$.
  In particular, the index of $\mathcal{O}$ in the maximal order is coprime to~$p$.

  Conversely, if $\End(E)$ is imaginary quadratic with $p$ split, say $(p) = \mathfrak{p}\cdot \mathfrak{p}'$, then the isogenies given by the kernels $E[\mathfrak p]$ and $E[\mathfrak p']$ have the same codomain isomorphic to $E$.
\end{proof}

In the case of a supersingular elliptic curve, here is an example later discussed in Section~\ref{ss_13_5_ex_subsec}:
Take the elliptic curve $E$ defined over $\FF_{13}$ with $j$-invariant equal to $5$.
This is the unique supersingular $j$-invariant.
This curve admits three necklaces defined over $\FF_{13}$, however the $j$-invariants of $E/C$ are all equal to $5$ as well.

In terms of isogeny volcanoes~\cite{sutherland_volcanoes}, this means that if the $j$-invariants are not distinct, the curve sits on the rim of the volcano with at least one pair of vertices with multiple connected edges.

\begin{prop}\label{distinct_j_prop}
  Let $E$ be an elliptic curve defined over a number field~$K$ and suppose $j(E)\not\in\{0,1728\}$.
  If the image of the representation $\rho\colon G_K\rightarrow \PGL\bigl(E[p]\bigr)$ is equal to a full nonoca, then the $j$-invariants $j(E/C)$ for $C$ cyclic in $E[p]$ are pairwise distinct.
\end{prop}

\begin{proof}
  Suppose that there exists two distinct subgroups $C$ and $C'$ of order $p$ such that $j(E/C)=j(E/C')$.
  Then, by the previous lemma, $E$ has complex multiplication by an imaginary quadratic order $\mathcal O$ in which $p$ splits.
  By Lemma~\ref{cm_classic_lem}, this implies that $\rho(G_K)$ is contained in the normaliser of a split Cartan subgroup.
  This is not possible since  $\rho(G_K)$ is a nonoca.
\end{proof}

\begin{prop}\label{distinct_j_overQ_prop}
  Suppose $p>3$.
Let $E$ be an elliptic curve defined over $\QQ$ with $j(E)\not\in\{0,1728\}$.
  If $E$ admits a necklace $\vv$ defined over $\QQ$, then the $j$-invariants $j(E/C)$ for $C$ cyclic in $E[p]$ are pairwise distinct.
\end{prop}

\begin{proof}
  As above, $j(E/C)=j(E/C')$ for $C\neq C'$ would imply that $E$ has complex multiplication by an imaginary quadratic field $F$ in which $p$ splits.
  But Lemma~\ref{cm_lem} contradicts then the existence of a necklace defined over~$\QQ$.
\end{proof}

\subsection{Curves with extra automorphisms}\label{extra_aut_subsec}
We discuss now the special cases when $\Aut(E)$ is strictly larger than~$\{\pm 1\}$.
Let $p>3$ be a prime number and $E$ be an elliptic curve defined over a field $k$ of characteristic different from~$p$, $2$, and $3$ and suppose that $j(E)=0$ or~$1728$.
Denote $F=\QQ(\sqrt D)$ and
\begin{align*}
  D&=-3,\quad & \mathcal{O}&=\mathcal{O}_F =\ZZ[\zeta],\quad & n&=3 & \text{ if } j(E)&=0\\
  D&=-4, & \mathcal{O}&=\mathcal{O}_F =\ZZ[i], & n&=2 & \text{ if } j(E)&=1728
\end{align*}
where $\zeta$ is a primitive cube root of unity in $\bar k$, and consider
\begin{equation*}
  E_{-3}\colon\quad y^2 = x^3 + 1\qquad\text{and}\qquad
  E_{-4} \colon \quad y^2 = x^3 + x.
\end{equation*}
As $E$ has complex multiplication by $\mathcal{O}$, it is a twist of~$E_D$.
It has an equation of the form $y^2=x^3+d$ if $j(E)=0$ or of the form $y^2=x^3+dx$ if $j(E)=1728$, with $d$ a $(2n)$-th power-free integer.

As $[-1]$ acts trivially on $\PP\bigl(E[p]\bigr)$,
the cyclic group $\Aut(E)/\{\pm 1\}$ of order $n$ acts on the set of necklaces on $E[p]$.
Let $\alpha=[\zeta]$ if $j=0$ and $\alpha=[i]$ if $j=1728$.
It induces an element $u\in\PGL\bigl(E[p]\bigr)$ which is of order~$n$.

\begin{prop}\label{autom_prop}
  Suppose $k=\QQ$ and $p>5$.
  If $j(E)=0$ assume that~$p>7$.
  Then there is a necklace in the $p$-torsion of $E$ defined over $\QQ$ if and only if $p\equiv 2\pmod{3}$ when $j(E)=0$ or $p\equiv 3\pmod 4$ if $j=1728$.
  This necklace $\vv^*$ is unique and the element $u$ acts as a rotation on it.
\end{prop}

\begin{proof}
  As $p>3$, the prime cannot be ramified in $\mathcal{O}$.
  As seen in the proof of Lemma~\ref{cm_classic_lem}, $\rho(G_F)$ is a subgroup of a Cartan subgroup $C_{\mathcal{O}},$ which is split if $p$ is split in $\mathcal O$ and non-split if $p$ is inert.
  From Stevenhagen's theorem cited in~\cite{bourdon_clark}, we know that the index of $\rho(G_F)$ in $C_{\mathcal{O}}$ is either $1$ or~$n$.

  Suppose $p$ is inert.
  Then $\rho(G_F)$ has at least $(p+1)/n$ elements.
  The lower bounds on $p$ ensure that this is greater than~$2$.
  As $\rho(G_F)$ contains more than $2$ elements, $\rho(G_\QQ)$ cannot be contained in a nonoca other than the normaliser of $C_\mathcal O$.
  Therefore there is a unique necklace defined over~$\QQ$.

  If $p$ is split, then there are at least $(p-1)/n>2$ elements.
  This implies that $\rho(G_\QQ)$ cannot be contained in a nonoca and hence there is no $\QQ$-rational necklace on~$E$.

The congruence conditions come from the fact that $p$ is inert in $\mathcal O$ if and only if $ \bigl(\tfrac{D}{p}\bigr)=-1$.

  If $E$ admits a necklace $\vv$, the element $u$ belongs to $\Aut_{\mathcal{O}}\bigl(E[p]\bigr) = \bigl(\mathcal O/p\mathcal O\bigr)^{\times}$ inside $\PGL\bigl(E[p]\bigr)$ and not just its normaliser.
  It then acts as a rotation on~$\mathfrak v$.
\end{proof}

One can find some curves $E/\QQ$ with $j(E)=0$ with a single and some that have two $\QQ$-rational necklaces for $p=5$.
Moreover, some have a $\QQ$-rational necklace for $p=7$ despite $7\equiv 1 \pmod{3}$.
Compare with Example~4.4 in~\cite{bourdon_clark}.

\begin{lem}\label{extraaut_lem}
  If $u$ is of order $2$, then its action on $\PP\bigl(E[p]\bigr)$ commutes with the action of the Galois group of $k$.
  If $u$ is of order~$3$ and $k$ contains the third roots of unity, then $u$ commutes with the Galois group.
  Otherwise the Galois group may invert $u$.
\end{lem}

\begin{proof}
  As $[-1]$ acts trivially on $\PP\bigl(E[p]\bigr)$, it is the group $\Aut_{\bar{k}}(E)/\pm 1$ that acts.
  As a Galois module this is either isomorphic to $\mu_4/\mu_2\cong \mu_2$ which has a trivial action by the Galois group or $\mu_6/\mu_2\cong\mu_3$ which is trivial only if $k$ contains the third roots of unity.
\end{proof}

Let us now consider the points $x \in X$ in $\pi^{-1}(\{0,1728\})$ for $\pi \colon X\longrightarrow \PP^1$.
Such a point can be represented by a pair $(E,\vv)$ with $E=E_{-3}$ or $E_{-4}$.
We have $x\in X(k)$ if and only if for each $\sigma \in G_k$, there exists an automorphism $\psi_\sigma\in\Aut(E)$ such that $\psi_\sigma(\vv)=\sigma(\vv)$.

Denote by $\Omega$ the $\Aut(E)/\{\pm 1\}$-orbit of $\vv$, which is either a sigleton or it contains $n$ necklaces.
The ramification index of $x$ above $\pi(x)$ is $\lvert\Aut(E)\rvert /\lvert \Stab(\vv)\rvert = \# \Omega$.
If $\Omega=\{\vv\}$, then the point $x$ is an elliptic point of~$X$ as it is unramified above $\pi(x)$.
We already counted these points in Proposition~12 in~\cite{necklaces}.
In this case, $x\in X(k)$ if and only if $\vv$ is defined over~$k$.

If $\# \Omega=n$, then $x$ is a ramified point of index $n$; it is represented by $(E,\mathfrak w)$ for each~$\mathfrak w\in \Omega$.
Such a point is in $X(k)$ if either each necklace in $\Omega$ is defined over~$k$, in which case $G$ is contained in the intersection of the corresponding $n$ nonocas, or $\Omega$ forms a single $G_{k}$-orbit.

\begin{itemize}
  \item Above $j=0$ there is no elliptic point if $p\equiv 1\pmod 3$ and only one if $p\equiv 2\pmod 3$.
  In this last case, the unique necklace $\vv^*$ is fixed by $u$, that is to say $u$ acts as a rotation of angle $\pm 2\pi/3$.
  It can be visualized by folding the necklace three times over itself.
  Since $u(C_{i})=C_{i+(p+1)/3}$, we have $j(E/C_{i})=j\bigl(E/C_{i+(p+1)/3}\bigr)$ for all~$i$.
  The other points are ramified points $x=\bigl[(E,\vv)\bigr]=\bigl[(E,u(\vv))\bigr]=\bigl[(E,u^2(\vv)\bigr)]$.
  \item Above $j=1728$, the automorphism $u$ may either fix or flip the necklace for an elliptic point.
  If $p\equiv 1\pmod 4$, there are no such necklaces fixed by~$u$, while for $p\equiv 3\pmod 4$ there is exactly one such necklace, as explained in Section~3.5.1 in~\cite{necklaces} where it appears as a pair of two oriented necklaces.
  On this necklace $u$ acts as a rotation of angle $\pi$ and we may picture it like folding this necklace $\vv^*$ twice on itself.
  Here $j(E/C_{i})=j\bigl(E/C_{i+(p+1)/2}\bigr)$ for all~$i$.

  Moreover, there are $\bigl(p-(\tfrac{-1}{p})\bigr)/2$ elliptic points such that $u$ acts as a reflection on the necklace.
  These are the flipped necklaces in Section~3.5.2 of~\cite{necklaces}.
  The number of them is obtained by subtracting from the number $e_2^+$ in Proposition~12 in~\cite{necklaces} the number of fixed necklaces.

  Finally, the other points in the fibre are ramified points $x=\bigl[(E,\vv)\bigr]=\bigl[(E,u(\vv))\bigr]$.
\end{itemize}

In the case where $k=\QQ$, the necklace $\vv^*$ above is the unique necklace of Proposition~\ref{autom_prop}.

\begin{lem}\label{autom_nb_lem}
  Let $p>7$.
  When $p\equiv 2\pmod 3$ the only point of $X(\QQ)$ above $j=0$ is $[(E_{-3},\vv^*)]$, and there is no such point if $p\equiv 1\pmod 3$.
  Similarly, $[(E_{-4},\vv^*)]$ is the unique point above $j=1728$ when $p\equiv 3\pmod 4$ and there is none when $p\equiv 1\pmod 4$.
\end{lem}

\begin{proof}
  As we are no longer concerned with the specific curve, but only the points in $X$, we may take the curve to be $E_D$.
  Recall that we denote $\pi \colon X\longrightarrow \PP^1$.
  An unramified point in $ \pi^{-1}(\{0,1728\})$ is defined over $\QQ$ if and only if it is represented by $(E_D,\vv)$ with $\vv$ defined over~$\QQ$:
  There is only one such point obtained for $\vv=\vv^*$ if $p\equiv -1\pmod D$, and none otherwise by Proposition~\ref{autom_prop}.

  We are left to prove that no ramified point $x=[(E_D,\vv)]$ is defined over~$\QQ$ when $\vv\neq\vv^*$, in other words we want to exclude the possibility that the $\Aut(E)/\{\pm 1\}$-orbit $\Omega$ of $\vv$ is a $G_\QQ$-orbit.
  For our particular curve $E_D$, the reduction is good at $p$ and $p$ is not ramified in $\mathcal O$.
  Therefore, Lemma~\ref{cm_classic_lem} implies that $G = \rho(G_{\QQ})$ is equal to the normaliser of a Cartan subgroup in $\PGL\bigl(E[p]\bigr)$.
  It is dihedral of order~$2(p+1)$ or order $2(p-1)$.
  Since $\vv$ is not fixed by $G$, its stabiliser under the Galois action is of order $2$ or $4$, since it is in the intersection of $G$ with the nonoca stabilizing~$\vv$.
  Hence the Galois orbit of $\vv$ has order $(p\pm 1)/2$ or $(p\pm 1)$.
  For $p>7$, the Galois orbit has at least $(7-1)/2 = 4$ elements, which implies that $\Omega$ cannot be such an orbit as $\Omega$ has at most $3$~elements.
\end{proof}

The fibres for smaller primes $p$ can be exceptionally.
For $p=7$, the fibre above $j=0$ contains one rational point and it is ramified despite $p\equiv 1 \pmod{3}$.
For $p=5$, the fibre above $j=0$ has two rational points, one ramified and one unramified.
This can be checked on the explicit models or with the same method as in the proof above.

\section{Representation of necklaces over number fields}\label{algo1_sec}
Let $E$ be an elliptic curve defined over a number field $K$, such that $j(E)\not\in\{0,1728\}$.
We suppose that $E$ admits a necklace $\vv$ defined over $K.$
The aim of the following algorithm is to represent $\vv$.
In this section, we suppose that the $j$-invariants $j(E/C)$ for $C$ cyclic in $E[p]$ are pairwise distinct (see Lemma~\ref{distinct_j_lem}): in this case, we may represent $\vv=(C_0,C_1,\dots,C_p)$ by an ordered list of the $j$-invariants
\[
  \bigl(j(E/C_0),\,j(E/C_1),\,\dots,\,j(E/C_p)\bigr).
\]
Actually, for algorithmic purposes, we will make a stronger assumption: We suppose that the image $G$ of the representation $\rho\colon G_K \to \PGL\bigl(E[p]\bigr)$ is the full nonoca stabilising~$\vv$.
In this case, the $j$-invariants are automatically distinct by Proposition~\ref{distinct_j_prop}.
See also Proposition~\ref{distinct_j_overQ_prop}.

Let $\mathfrak q\mid q$ be a prime ideal of $K$ of good reduction for $E$ with~$q\neq p$.
Denote by $\FF_{\mathfrak q}$ the residue field of $K$ at $\mathfrak q$ and $\widetilde E$ the reduction of $E$ modulo~$\mathfrak q$.
Then $E[p]\cong \widetilde E[p]$ as $G_K$-modules, where $G_K$ acts on $\widetilde E[p]$ via the canonical surjection $G_K\to G_{\FF_{\mathfrak q}}$.
On $\PP\bigl(E[p]\bigr)\cong\PP\bigl(\widetilde E[p]\bigr)$ this action will factor through the \emph{$p$-isogeny field}~$L$, that is to say the smallest extension of~$K$ over which all cyclic subgroups of $E$ of order $p$ are defined.
Since $E$ has good reduction at~$\mathfrak q$, the inertia group at~$\mathfrak q$ acts trivially on $E[p]$.
Hence the action of $G_K$ on $\widetilde E[p]$ is cyclic generated by any choice of a Frobenius element $\Frob_{\mathfrak Q}\in\Gal(L/K)$ for $\mathfrak Q\mid \mathfrak q$ an ideal of $L$.

The basic idea of the algorithm described below as Algorithm~\ref{nf_alg} is the following:
\begin{description}
  \item{Step 1}:
  We calculate $L$ as the splitting field of the polynomial $f(x) = \Phi_p\bigl(j(E),x\bigr)$, where $\Phi_p$ is the standard modular polynomial for~$Y_0(p)$.
  These polynomials have been calculated by Bröker, Lauter, and Sutherland as explained in~\cite{broker_lauter_sutherland} and can simply be read off a file which can be found on Sutherland's webpage.
  \item{Step 2}: We determine a prime ideal $\mathfrak q$ in $K$ such that the roots of $f$ have distinct reduction modulo a prime $\mathfrak{Q}$ in $L$ above $\mathfrak{q}$ and such that a Frobenius $\Frob_{\mathfrak Q}\in\Gal(L/K)$ is the class of an element of order~$p+1$.
  More precisely, we want $\Frob_{\mathfrak Q}$ to have characteristic polynomial on $E[p]$ equal to the minimal polynomial of~$\gamma$.
  To obtain this, we try the primes $\mathfrak{q}$ of good reduction such that $\Norm(\mathfrak q) \equiv \Norm(\gamma) \pmod{p}$ and $a_{\mathfrak{q}}\equiv \Tr(\gamma)\pmod{p}$, each time checking if the reduced roots of $f$ are distinct.
  Here $a_{\mathfrak{q}} = \Norm(\mathfrak{q}) +1 - \#\tilde E(\FF_{\mathfrak{q}})$.
  \item{Step 3}: In the last step, we order the roots $j(E/C)\in L$ of $f$ according to the necklace $\vv$:
  We pick a first root $j_0$ among them.
  Then we pick $j_1$ to be the unique root whose reduction modulo $\mathfrak{q}$ is $\Frob_{\mathfrak{Q}}(j_0)$ in the residue field.
  Then $j_2$ and so forth.
\end{description}
\begin{algorithm}\label{nf_alg}
  \SetKw{KwFrom}{from}
  \SetProgSty{}
  \DontPrintSemicolon
  \KwIn{$E/K$ elliptic curve as above, a prime $p$, and a generator $\gamma$ of $\FF_{p^2}^{\times}$.}
  \KwOut{A list $(j_0,j_1,\dots,j_p)$ of elements in a number field representing the necklace $\vv$ in $E[p]$}
  \;
  Read $\Phi_p$ in the file described above\tcc*[f]{Step 1}\;
  $f(x) \gets \Phi_p(j(E), x)\in K[x]$\;
  $L\gets $ the splitting field of $f$\;
  $J\gets $ the set of all roots of $f$ in $L$\;
  Set $t\gets \operatorname{Tr}(\gamma)$ and $n\gets \operatorname N(\gamma)$\tcc*[f]{Step 2}\;
  \Repeat{the reduction of elements in $J$ are distinct modulo $\mathfrak{Q}$}{
    \Repeat{ $\operatorname{N}(\mathfrak q) \equiv n \pmod{p}$ and $a_{\mathfrak q} \equiv t \pmod{p}$}{ Advance to the next prime ideal $\mathfrak q\nmid p$ in $K$ for which $E$ has good reduction}
    Pick a prime $\mathfrak{Q}$ in $L$ above $\mathfrak{q}$\;}
  $j_0 \gets$ one element in $J$\tcc*[f]{Step 3}\;
  \For{$k$ \KwFrom $1$ \KwTo $p$}{
    $y\gets \Frob_{\mathfrak{Q}}(j_{k-1}+\mathfrak{Q})$\;
    Set $j_k$ to be the element in $J$ that reduces to $y$ modulo $\mathfrak{Q}$.}
  \Return the necklace $(j_0,j_1,\dots, j_p)$
\caption{Computing necklaces over number fields}
\end{algorithm}

The bottleneck of the algorithm is the complete factorisation of $f$ in the field $L$.

\begin{rmk}
  The condition on $\Frob_\mathfrak Q$ to be of order $p+1$ implies that the prime of good reduction $\mathfrak q$ for $E$ is ordinary, and by Lemma~\ref{distinct_j_lem}, the condition on the roots of $f\pmod {\mathfrak Q}$ implies that $p$ is not split in $\End(\tilde E)$.
\end{rmk}

\section{Reduction of a necklace}\label{red_sec}
Let $E$ be an elliptic curve defined over a number field $K$ with $j(E)\not\in\{0,1728\}$.
Suppose $\vv$ is a necklace in $E[p]$ defined over $K$.

\subsection{Good reduction}
Let $\lambda$ be a prime ideal in $K$ not dividing $p$, $\FF_\lambda$ the residue field at $\lambda$ and suppose that $E$ has good reduction at~$\lambda$.
Then the point $[(E,\vv)]\in Y(K)$ can be reduced to a point $\bigl[(\tilde E, \tilde {\vv})\bigr]\in Y(\FF_{\lambda})$.
If the necklace is given by the above Algorithm~\ref{nf_alg}, then we can try obtaining a representation of $\tilde{\vv}$ simply as follows:
Pick a prime $\mathfrak{L}$ above $\lambda$ in the $p$-isogeny field $L$ and reduce the values $j_k = j(E/C_k)$ modulo~$\mathfrak{L}$.
If we are lucky the reduced values are distinct.
We would then have a list $(\tilde{j}_0,\tilde{j}_1,\dots, \tilde{j}_p)$ belonging to a finite extension of $\FF_{\lambda}$ representing the necklace $\tilde{\vv}$ by listing the $j$-invariants of the curves which are the codomains of $p$-isogenies leaving the reduced curve~$\tilde{E}/\mathbb{F}_{\lambda}$.

However, we may be unlucky.
If the reduction of $E$ at $\lambda$ is supersingular, we may expect that the $j$-invariants do not reduce to distinct elements modulo $\mathfrak{L}$.
Lemma~\ref{distinct_j_lem} explains that even when the reduction at $\lambda$ is good ordinary, the $j$-invariant may no longer pairwise distinct in the reduction.

In this situation, we need to do extra work.
First, we can write down explicitly isogenies $\varphi_k\colon E \to E/C_k$ given that we know the degree and the two elliptic curves involved.
There are two choices for the isogeny, $\varphi_k$ and $-\varphi_k$, with the same kernel $C_k$ as we have no extra automorphisms by assumption.
However, this choice will not matter as we care for $C_k$ rather than for~$\varphi_k$.
We obtain this way a kernel polynomial $f_k$ defining the cyclic subgroup $C_k$ as a subgroup scheme of $E$ and this does not depend on the above choice of~$\varphi_k$.
This kernel polynomial can be reduced modulo~$\mathfrak{L}$.
We obtain a list of polynomials $\bigl(\tilde{f}_0, \tilde{f}_1,\dots, \tilde{f}_p\bigr)$ defined over a finite field.
This represents the reduced necklace~$\tilde{\vv}$.
We could also reduce the isogenies and represent it as a list $\bigl( \tilde{\varphi}_0,\tilde{\varphi}_1,\dots,\tilde{\varphi}_p\bigr)$.

\begin{algorithm}\label{red_alg}
  \SetProgSty{}
  \SetKw{KwFrom}{from}
  \DontPrintSemicolon
  \KwIn{$E/K$ elliptic curve as above, a prime $p$, necklace represented by distinct $j$-invariants $(j_0,j_1,\dots,j_p)$ and a prime $\lambda$.}
  \KwOut{Either a list of $p+1$ distinct elements $(\tilde{j}_0, \dots, \tilde{j}_p)$ or an ordered list of $p+1$ polynomials $\bigl(\tilde{f}_0, \tilde{f}_1,\dots, \tilde{f}_p\bigr)$ defined over a finite extension of $\FF_{\lambda}$}
  \;
  $L\gets$ $p$-isogeny field of $E$\;
  Pick a prime $\mathfrak{L}$ in $L$ above $\lambda$\;
  Reduce the $j$-invariants modulo $\mathfrak{L}$ to $\tilde J = (\tilde{j}_0, \tilde{j}_1,\dots, \tilde{j}_p)$\;
  \If{all reduced values in $\tilde{J}$ are distinct}{
    \Return $\tilde{J}$}
  \For{$j$ \KwFrom $0$ \KwTo $p$}{
      Determine the isogeny $\varphi_k\colon E \to E/C_k$ given the codomain and degree\;
      Calculate the kernel polynomial $f_k$ defining $C_k$ in $L[x]$\;
      Reduce the kernel polynomial modulo $\mathfrak{L}$\;
    }
    \Return $\bigl(\tilde f_0, \tilde f_1,\dots,\tilde f_p\bigr)$\;
  \caption{Reduce a necklace}
\end{algorithm}

\subsection{Bad reduction}
Suppose now that $E$ has bad reduction at $\lambda$.

First, if the reduction of $E$ is potentially multiplicative, i.e., $j(E)$ has negative valuation at $\lambda$, then the reduction will be one of the cusps of $X$ over $\FF_{\lambda}$.
We could use the description of necklaces on Tate curves as in the proof of Proposition~6 in~\cite{necklaces} to decide which of the $\tfrac{p-1}{2}$~cusps the point reduces to, but we have not implemented this.

Otherwise, $j(E)$ is integral at $\lambda$.
We will need to reduce the $j$-invariants $j(E/C)$, which belong to the $p$-isogeny field $L$, or the kernel polynomials, whose coefficients are in~$L$.
Over $L$ the elliptic curve $E$ will acquire good reduction at primes above $\lambda$ due to the following lemma.
Therefore the $j$-invariants and polynomials can be reduced at the chosen prime $\mathfrak{L}$ above~$\lambda$.
The obtained reduced elliptic curve $\tilde E/\FF_{\mathfrak{L}}$ will admit a non-singular model defined over $\FF_{\lambda}$ as its $j$-invariant lies in $\FF_{\lambda}$.
Thus the reduced information will look exactly like in the case of good reduction.

\begin{lem}
  Let $E$ be an elliptic curve defined over a local field $k$ of residual characteristic~$\ell$.
  If $E$ admits an isogeny $E\to E'$ defined over $k$ of prime degree $p>3$ with $(p,\lambda)=1$, then $E$ has semistable reduction.
\end{lem}

\begin{proof}
  Let $C$ be the kernel of the isogeny viewed as a subgroup scheme of the Néron model~$\mathcal{E}$.
  Assume that $\mathcal{E}$ has additive reduction.
  As the group of components has order at most $4$ and $p>3$, the subgroup $C$ lies in the connected component of the identity~$\mathcal{E}^0$.
  Since $p$ is coprime to $\ell$, the special fibre of $C$ is étale.
  This is impossible as $\mathbb{G}_a$ over a field of characteristic $\ell$ has no subgroup of order~$p$.
\end{proof}

\subsection{Reduction of necklaces on CM elliptic curves}\label{cm_alg_subsec}
We now present a quicker algorithm to calculate reductions of necklaces on \textsc{cm}~elliptic curves avoiding Algorithm~\ref{nf_alg}.
Instead we use the reduction of the endomorphism ring.

Let $E$ be an elliptic curve defined over a number field~$K$.
Assume that $E$ has complex multiplication by an order $\mathcal{O}$ of an imaginary quadratic field $F$ and that $p$ is inert in~$\mathcal O$.
Then, as discussed in Section~\ref{cm_subsec}, there exists a special necklace $\vv^*$ defined over $K$ on $E[p]$ coming from the fact that $E[p]$ is a free $\mathcal{O}/p\mathcal{O}$-module of rank~$1$.
This is the necklace fixed by the non-split Cartan subgroup $C_\mathcal O$ in $\PGL\bigl(E[p]\bigr)$ coming from the $\mathcal{O}$-structure.
Although we do not elaborate further in this article, we note that the point $[(E,\vv^*)]\in X(K)$ is a Heegner point.

If $\rho(G_K)=C_{\mathcal O}$, then $\vv^*$ is the unique necklace on $E[p]$ defined over~$K$.
This is the case for instance if $K=\QQ$ by Lemma~\ref{cm_lem} or more generally if $p\nmid f\Delta_{FK}$ is a place of good reduction for $E$ by Lemma~\ref{cm_classic_lem}.

Let $\lambda\nmid p$ be a prime ideal of $K$ of good reduction for~$E$.
Denote by $\FF_\lambda$ the residue field of $K$ at~$\lambda$.
By reducing $\vv^*$ modulo $\lambda$, we obtain the necklace~$\tilde\vv^*$ on the reduction~$\tilde E$.
The aim of this section is to present an algorithm to calculate $\tilde\vv^*$ from only data of the \textsc{cm}~elliptic curve $E$ without having to calculate $\vv^*$ beforehand, that is without executing the costly Algorithm~\ref{nf_alg}.
By Proposition~3.4 in~\cite{rubin}, the reduction map $\mathcal O\cong\End(E) \hookrightarrow\End(\tilde E)$ is injective.
Hence $\tilde\vv^*$ is the unique necklace on $\tilde E[p]$ stabilised by $ (\mathcal O/p\mathcal O)^\times$.
It is then given by the action of any element of $\mathcal O\hookrightarrow\End(\widetilde E)$ which modulo $p$ maps to the chosen generator $\gamma$ of $\FF_{p^2}^\times$ through the isomorphism with~$(\mathcal O/p\mathcal O)^\times$.

The method is presented in Algorithm~\ref{cm_alg}, for elliptic curves over $\QQ$.

\begin{algorithm}\label{cm_alg}
  \LinesNumbered
  \SetProgSty{}
  \SetKw{KwFalse}{False}
  \SetKw{KwTrue}{True}
  \SetKw{Kwin}{in}
  \DontPrintSemicolon
  \KwIn{Two distinct primes $p$ and $\ell$, an elliptic curve $E/\QQ$ with complex multiplication, and a unique necklace $\vv$ in its $p$-torsion}
  \KwOut{A representation of the necklace $\tilde{\vv}$ on $\tilde E$ over $\FF_{\ell}$ }
  \vspace{1ex}
  Determine $\End(E)$ and find a $\ZZ$-basis $\{1,\psi\}$ for it\;
  Find an element $\varphi =a+ b\, \psi\in\End(E)$ such that the reduction in $\End(E)/p\End(E) \cong \FF_{p^2}$ maps to $\gamma$\;
  Reduce $\varphi$ to $\tilde\varphi \in \End(\tilde{E})$ defined over $\FF_{\ell^2}$\;
  Determine the $p$-isogeny field for $\tilde E$ and all $p$-isogenies leaving $\tilde E$\;
  Order them as in the necklace $\tilde{\vv}$ by acting with $\tilde\varphi$ on $p$-isogenies\;
  \caption{Construct the reduced necklace for a curve with complex multiplication}
\end{algorithm}

Again some remarks should be made.
First of all, we opt to represent isogenies, and in particular endomorphisms on elliptic curves, as formal sums of compositions of easier isogenies.
This is done effectively in Sage~\cite{sagemath}.
For instance only $[a]$, $[b]$ and $\psi$ are used and the actual rational map $\varphi=a+b\psi$ is never directly calculated.

For line~1, one needs to construct an endomorphism $\psi$ which is not in $\ZZ$.
For $j=0$ and $j=1728$ one can take an automorphism of $E$ other than $[\pm 1]$; this is a simple change of variables in the Weierstrass equation.
For all other curves, we pick a small prime $q\neq \ell$ which splits in~$\mathcal{O}$.
There is an endomorphism of degree $q$ on $E/F$ which can be constructed explicitly.
The hardest case for this is for the curve with $F=\QQ\bigl(\sqrt{-163})$ and $\ell=41$, but the requested endomorphism of degree~$q=43$ is not difficult to calculate either.
The identification $[\cdot]\colon\mathcal{O} \to \End(E)$ obtained in this way has to be normalised such that $[a]^*(\omega) = a\cdot \omega$ for a differential $\omega$ on~$E$.

In line~2, any lift of our $\gamma$ to $\mathcal{O}$ will do.
However, the reduction to an $\FF_p$-scalar multiple of $\gamma$ will have the same action on the pearls, and hence any lift of those will also work.
The possible elements in $\End(E)$ are all elements in a subgroup of index $p$ which do not belong to~$p\End(E)$.
Therefore the values of $a$ and $b$ can be chosen fairly small compared to~$p$; although in practice these values do not matter too much as we will work with the formal sum as explained above.

In step~3, we reduce the endomorphism modulo a prime ideal $\lambda$ in $F$ above~$\ell$.
In practice, as the isogeny is given as a sum of compositions, it is best to reduce these components and still represent it as a sum of composition over the residue field~$\FF_\lambda$.
As we have avoided that $\psi$ has degree divisible by $\ell$, the reduced isogeny is obtained as a composition of separable isogenies.

In the actual implementation in~\cite{implementation} of step~4, we construct the $p$-division field $\FF_{\ell}\bigl(\tilde{E}[p]\bigr)$ and consider the action of $\tilde{\varphi}$ on the $p$-torsion points.
This is because it would be extra work to implement the Galois action directly on isogenies.
Since points and isogenies are efficiently implemented over finite fields, this does not significantly reduce the speed of this algorithm.

We have assumed so far that the reduction of $E$ is good at~$\ell$.
However, it is easy to pass by.
From Corollary~5.22 in~\cite{rubin}, we know that there is an elliptic curve $E'$ over $F$ which is isomorphic to $E$ and which has good reduction at the fixed prime $\ell$.
In practice, the curve $E'$ can be obtained by a quadratic twist.

\section{Comparing reduced points on the modular curve}\label{comparing_sec}

\subsection{Algorithm for comparing points over a finite field}
In the previous section, we saw how we can obtain points in $Y$ defined over the residue field $\FF_{\lambda}$ of a number field.
We will present an algorithm to test if two points $x_1=\bigl[(E_1,\vv_1)\bigr]$ and $x_2=\bigl[(E_2,\vv_2)\bigr]$ in $Y(K)$ reduce to the same point in $Y(\FF_{\lambda})$ at a given prime ideal $\lambda$ in $K$ not dividing~$p$.
For this we will explain Algorithm~\ref{red_equal_alg} which compares points on $Y$ over finite fields.

Note first that, as seen in Lemma~\ref{number_of_necklaces_lem}, in some cases there is a single necklace on the elliptic curve, in which case we only need to check if the curves are isomorphic.
But most elliptic curves over a finite field will have more than one necklace defined over that field.

We are given two elliptic curves $E_1$ and $E_2$ over a finite field $\FF$ of characteristic $\ell\neq p$, respectively endowed with a necklace $\vv_1$ and $\vv_2$ also defined over~$\FF$.
In practice, these necklaces are either given as a list of distinct $j$-invariants in an extension of $\FF$ or as a list of kernel polynomials with coefficients in an extension of~$\FF$.
We treat here first the case when both are given by list of distinct $j$-invariants.

First, we check that $E_1$ and $E_2$ are $\bar{\FF}$-isomorphic by checking if $j(E_1) = j(E_2)$ in~$\FF$.
Next, we can check if the $p$-isogeny fields of $E_1$ and $E_2$ are isomorphic.
Finally, we may use an isomorphism between them to check if the list of $j$-invariants represents the same necklace by checking whether one list is a cyclic shift or a cyclic shift composed with reversing of the other list.

\begin{algorithm}\label{red_equal_alg}
  \SetProgSty{}
  \SetKw{KwFalse}{False}
  \SetKw{KwTrue}{True}
  \SetKw{Kwin}{in}
  \DontPrintSemicolon
  \KwIn{Two elliptic curves $E_1$ and $E_2$ over a finite field $\FF$ each with a necklace $\vv_1$ and $\vv_2$ given by a list $J_1$ and $J_2$ of $p+1$ elements in an extension of $\FF$.}
  \KwOut{Boolean deciding if $(E_1,\vv_1)$ and $(E_1,\vv_2)$ represent the same point in $Y(\FF)$}
  \;
  \If{$j(E_1)\neq j(E_2)$}{\Return \KwFalse}
  $F_1\gets$ the $p$-isogeny field of $E_1$\;
  $F_2\gets$ the $p$-isogeny field of $E_2$\;
  \If{$F_1 \not\cong F_2$}{\Return \KwFalse}
  Identify $F_1$ and $F_2$ and use it to convert elements in $J_1$ and $J_2$ to the same field\;
  \If{$J_1$ differs from $J_2$ as a set}{\Return \KwFalse}
  \If{the permutation from $J_1$ to $J_2$ is either a $(p+1)$-cycle or a $(p+1)$-cycle composed with reversing the order}{\Return \KwTrue}
  \Else{\Return \KwFalse}
  \caption{Test if two reduced necklaces are equal}
\end{algorithm}

If the $j$-invariants are not all distinct, we have to work with two lists of polynomials $(f_0, f_1,\dots, f_p)$ instead.
The basic comparison is as above: First checking if the curves are isomorphic, then if they have the same $p$-isogeny field and, finally, if the two kernel polynomial lists are linked by the correct permutation.
Note that one has to make sure that the polynomials are consistently normalised to compare correctly if they give the same subgroup under an isomorphism of $E_1$ and $E_2$.

In this last step, one has to treat the special case that the $j$-invariant may be $0$ or $1728$; a case that can never appear in the first version of the algorithm as the $j$-invariants $j(E/C_k)$ will not be distinct.
For this situation, one needs to account for extra automorphisms as explained in Section~\ref{extra_aut_subsec}.
In practice we check if the second list of polynomials is the right sort of permutation of $(u^*(f_0), u^*(f_1),\dots,u^*(f_p))$ for any automorphism $u$ of $E_1$ defined over the same field as~$f_k$.
Here $u^*(f_k)$ is the polynomial defining~$u(C_k)$.

\begin{rmk}
  We have fixed a generator $\gamma$ at the start and our comparison assumes that both necklaces were constructed with the same choice of~$\gamma$.
  If this were not the case, one could check equality by finding an $n$ coprime to $p+1$ and a $k$ such that $j(E_1/C_i) = j(E_2/C_{ni+k})$ for all $i$ and similar for the composition with reversing the order.
  See Lemma~1 in~\cite{necklaces}.
\end{rmk}

\subsection{Injectivity of reduction}\label{inj_subsec}
Let $p$ and $\ell$ be two distinct prime numbers larger than $3$.
Denote again by $X$ a model of $X$ over $\ZZ[\tfrac{1}{p}]$ and by $\red_\ell \colon X(\QQ) = X\bigl(\ZZ[\tfrac{1}{p}]\bigr) \to X(\FF_\ell)$ the reduction map.

For $p=5$ or $p=7$, when $X\cong \PP^1$ as a $\ZZ\bigl[\tfrac{1}{p}\bigr]$-scheme, this map is the obvious surjective reduction map.
Similar for $p=11$, when the genus is $1$ and $X(\QQ)$ is an elliptic curve with positive rank.
The situation is different for $p>11$, when $X(\QQ)$ is finite.

It is conjectured that for $p>11$, the set $X(\QQ)$ is equal to the set $\cmset$ of rational points represented by $(E,\vv)$ such that $E$ has complex multiplication.
In view of this conjecture, we are interested in~$\red_\ell\vert_{\cmset}$.
In particular, we can ask for which $\ell$ this map is injective.
In other words, we discuss the question for which $\ell$ are there two distinct points $x=[(E,\vv)]$ and $x'=[(E',\vv')]$ in $X(\QQ)$ both with complex multiplication and having the same reduction in $X(\FF_{\ell})$.

Suppose that $x$, $x'$, and $\ell$ are as above and~$p>7$.
Since, by Lemma~\ref{cm_lem} and Lemma~\ref{autom_nb_lem}, there is only one point in $\cmset$ with a given $j$-invariant, we represent $x=\bigl[(E,\vv)\bigr]$ and $x'=\bigl[(E',\vv')\bigr]$ with~$j(E)\neq j(E')$.
Write $\widetilde E$ for the common reduction modulo~$\ell$.
The prime number $\ell$ must divide $j(E)-j(E')$, hence $\ell$ is in the finite list $\mathcal L$ of all prime divisors of the finitely many differences of \textsc{cm} $j$-invariants over~$\QQ$:
\begin{multline*}
  \mathcal L= \bigl\{3\leq\ell \leq 127 \ : \ \ell\text{ is prime }\bigr\} \ \cup \\ \bigl\{137,139,157,163,173,193,197,211,229,233\bigr\} \ \cup \\
  \bigl\{241,257,277,283,293,317,331,389,433,571,643,997\bigr\}
\end{multline*}

\begin{lem}
  Let $(E,\vv)$ and $(E',\vv')$ be two \textsc{cm}~points in $X(\QQ)$ and $\ell\neq p$ an odd prime such that their common reduction $\tilde E$ at $\ell$ is ordinary.
  Then the reductions $(\tilde E,\tilde{\vv})$ and $(\tilde{E},\tilde{\vv}')$ in $X(\FF_{\ell})$ are distinct.
\end{lem}

\begin{proof}
  We can assume that $E=E_D$ and $E'=E_{D'}$ from our list in Table~\ref{cm_overQ_tab}.
  From this table, we see that no ordinary prime $\ell$ divides the conductor of $\End(E)$ as its divisors are ramified in $\End(E) \otimes \QQ$.
  By Deuring's result, as stated in Theorem 12 of section 13.4 in~\cite{lang}, it follows that $\End(E)\cong \End(\tilde E)\cong\End(E')$ and hence~$j(E)=j(E')$ since the endomorphism rings have class number $1$.
\end{proof}

Also we have checked algorithmically that the only case when $(E,\vv)$ and $(E',\vv')$ have equal reduction and both good ordinary reduction is when $\ell=2$ and $E$ and $E'$ are the curves $E_{-7}$ and~$E_{-28}$.

For each fixed $p$, as we will explain below, with our Algorithms~\ref{red_equal_alg} and~\ref{cm_alg}, we can effectively determine all the finitely many primes $\ell\in \mathcal L$ for which $\red_\ell\vert_{\cmset}$ is not injective.

\begin{prop}\label{inj_prop}
  Among the curves $X$ for primes $3< p <50$, there are only eight cases of $(\ell,E,\vv,E',\vv')$ such that $x=\bigl[(E,\vv)\bigr]$ and $x'=\bigl[(E',\vv')\bigr]$ are in $\cmset \subset X(\QQ)$ and such that $\red_\ell(x)=\red_\ell(x')$.
  They are all listed in the following table.
  \begin{center}
  \begin{tabular}{rrrrrl}
    \toprule
    $p$ & $\ell$ & $j$ & $\#X_j$ & $r_j$ & curves \\
    \midrule
    5 & 7 & 6 & 4 & 6 & $\bigl(E_{-7},E_{-163}\bigr)$; $\bigl(E_{-43},E_{-67}\bigr)$\\
    5 & 11 & 1 & 2 & 3 & $(E_{-27}, E_{-163})$\\
    5 & 13 & 5 & 4 & 5 & $(E_{-28}, E_{-67})$; $(E_{-8},E_{-163})$ \\
    5 & 17 & 8  & 4  & 4  & $(E_{-12}, E_{-163})$  \\
    7 & 5 & 0 & 3 & 4 &  $(E_{-8}, E_{-163})$ \\
    7 & 13 & 5 & 3 & 4 & $(E_{-67}, E_{-163})$ \\
    \bottomrule
  \end{tabular}
  \end{center}
\end{prop}

In this table, $X_j$ denotes the fibre in $X(\FF_{\ell})$ of $X \to \PP^1$ above~$j\in\FF_{\ell}$.
The column~$r_j$ counts the number of \textsc{cm}~points in $X(\QQ)$ which have this $j$-invariant modulo~$\ell$.
In all but one of the above cases, the number $r_j$ is larger than $\#X_j$ and hence the reduction cannot possibly be injective.
We have not found an example of non-injectivity when $r_j < \#X_j$, but also no explanation as to why this should not occur.

\begin{rmk}
  In all those cases, as seen in Lemma \ref{number_of_necklaces_lem}, since $\tilde E$ is supersingular, the image of the Frobenius has order $2$ in $\PGL\bigl(\tilde E[p]\bigr)$.
  If it is split, it flips $(p-1)/2$ necklaces, acting like a reflection of axis passing through two antipodal pearls.
  If it is non-split, it fixes one necklace as an angle $\pi$ rotation, and it flips $(p+1)/2$ necklaces acting like a reflection of axis passing between two couples of antipodal pairs.
  See Section~\ref{ex_sec} for an illustration of this.
\end{rmk}

The proof of the proposition is an explicit computer calculation using the implementation of Algorithm~\ref{cm_alg}.
Many of these instances are also explained in Section~\ref{ex_sec} as examples of necklaces.
For $p=5$, $7$, and $11$, we have also verified this result on the models as we can explain here in some details.
For $p=5$, in the \href{https://beta.lmfdb.org/ModularCurve/Q/5.10.0.a.1/}{\textsc{lmfdb}} there is an explicit description of the $j$-map from a model of $X$ as $\PP^1$ over $\ZZ\bigl[\frac{1}{30}\bigr]$, which can be characterised by saying that the nine \textsc{cm}~points have the following coordinates: $ E_{-3}= (-1:2)$, $E_{-7}=(1:0)$, $E_{-8} = (0:1)$, $E_{-12}=(1:2)$, $E_{-27}=(1:1)$, $E_{-28}=(-1:4)$, $E_{-43}= (-1,3)$, $E_{-67}=(5:6)$, and $E_{-163}=(-13:42)$.
From this is it easy to verify the assertion for $p=5$ made in the top four lines of the table.

For $p=7$, there is a model found in the \href{https://beta.lmfdb.org/ModularCurve/Q/7.21.0.a.1/}{\textsc{lmfdb}} which gives the \textsc{cm}~points as $E_{-4}= (0:1)$, $E_{-8}=(1:0)$, $E_{-11}=(-8:5)$, $E_{-16}=(-16:5)$, $E_{-43}=(8:5)$, $E_{-67}=(-4:5)$ and $E_{-163}=(-72:25)$.
The only congruences modulo prime $\ell>3$ are $(1:0)\equiv (-72:25)$ modulo $5$ and $(-4:5)\equiv (-72:25)\pmod{13}$ as expected.

Finally for $p=11$, the curve $X$ is isomorphic to the elliptic curve with Cremona label \href{https://www.lmfdb.org/EllipticCurve/Q/121b1/}{121b1} given by
\begin{equation*}\label{121b1_eq}
  y^2 + y = x^3 - x^2 -7\,x - 10
\end{equation*}
as found in~\cite[Prop. 4.3.8.1]{ligozat}.
This curve has $X(\QQ) = \ZZ\, Q$ with $Q=(4,5)$.
There is an isomorphism such that $E_{-12}$ maps to $O$.
The other \textsc{cm}~points map to $E_{-3}=\bigl(\tfrac{5}{4},\tfrac{7}{8}\bigr) = 3Q$,  $E_{-4}=(2,0)=2Q$, $E_{-16}=(4,-6)=-Q$,  $E_{-27}=(2,-1)=-2Q$,  $E_{-67}=(4,5)=Q$, and $E_{-167}=(-2,3)=4Q$.
The question if two points reduce to the same point modulo $\ell$ becomes the question if the coordinates of their difference have $\ell$ as a prime in the denominator.
The only points involved in such differences are the above points and $\pm 6Q$ and $\pm 5Q$.
Since $6Q = \bigl(\tfrac{25}{16},-\tfrac{85}{64}\bigr)$ and $5Q=\bigl(-\tfrac{8}{9},-\tfrac{118}{27}\bigr)$ are also $\{2,3\}$-integral, we can confirm that no two \textsc{cm}~points reduce to equal points modulo any prime $\ell > 3$.

Maybe it should not be surprising that we found very few instances of non-injectivity as the number of points of $X(\FF_{\ell})$ even for small $\ell$ increases quickly as $p$ grows.
With the table presented in Appendix~\ref{nu_pts_app}, it is easy to count the number of points in the reduction. 

In the following table, we list $\#X(\FF_{\ell})$ for $3<p\neq \ell < 50$.
The number $r=\#\cmset$ counts the number of \textsc{cm}~points in $X(\QQ)$.
The boldface number are the ones for which the reduction map on \textsc{cm}~points is not injective.

\newcommand\hl[1]{\textcolor{darkblue}{\mathbf{#1}}}
\begin{table}[H]
\centering
\begin{tabular}{cc|*{13}{r}}
  \toprule
  $p$ & $r$ & $\ell = 5$ & $7$ & $11$ & $13$ & $17$ & $19$ & $23$ & $29$ & $31$ & $37$ & $41$ & $43$ & $47$  \\
  \midrule
  $5$ & $9$ & & $\hl{8}$ & $\hl{12}$ & $\hl{14}$ & $\hl{18}$ & $20$ & $24$ & $30$ & $32$ & $38$ & $42$ & $44$ & $48$ \\
  $7$ & $7$ & $\hl{6}$ & & $12$ & $\hl{14}$ & $18$ & $20$ & $24$ & $30$ & $32$ & $38$ & $42$ & $44$ & $48$  \\
  $11$ & $7$ & $9$ & $8$ & & $14$ & $18$ & $20$ & $33$ & $30$ & $37$ & $31$ & $42$ & $44$ & $60$  \\
  $13$ & $7$ & $10$ & $11$ & $20$ & & $20$ & $24$ & $29$ & $31$ & $37$ & $26$ & $49$ & $31$ & $66$  \\
  $17$ & $7$ & $11$ & $15$ & $24$ & $13$ & & $23$ & $41$ & $45$ & $45$ & $27$ & $54$ & $37$ & $63$  \\
  $19$ & $7$ & $14$ & $13$ & $27$ & $14$ & $27$ & & $40$ & $45$ & $41$ & $38$ & $54$ & $39$ & $69$  \\
  $23$ & $7$ & $13$ & $16$ & $25$ & $22$ & $34$ & $30$ & & $47$ & $51$ & $31$ & $69$ & $44$ & $93$  \\
  $29$ & $8$ & $16$ & $18$ & $37$ & $20$ & $37$ & $46$ & $52$ & & $64$ & $30$ & $85$ & $57$ & $100$  \\
  $31$ & $8$ & $20$ & $18$ & $38$ & $21$ & $42$ & $42$ & $56$ & $65$ & & $32$ & $90$ & $51$ & $99$  \\
  $37$ & $4$ & $21$ & $21$ & $41$ & $23$ & $50$ & $49$ & $73$ & $72$ & $69$ & & $83$ & $65$ & $111$  \\
  $41$ & $8$ & $22$ & $28$ & $49$ & $27$ & $48$ & $54$ & $74$ & $75$ & $83$ & $35$ & & $59$ & $129$  \\
  $43$ & $4$ & $24$ & $23$ & $50$ & $27$ & $56$ & $49$ & $78$ & $75$ & $78$ & $40$ & $110$ & & $140$  \\
  $47$ & $8$ & $25$ & $30$ & $49$ & $34$ & $60$ & $56$ & $78$ & $79$ & $81$ & $48$ & $116$ & $69$ &  \\
  \bottomrule
\end{tabular}
\end{table}

Note that for the top two rows when $p=5$ or $p=7$, we have $\#X(\FF_{\ell}) = \ell+1$; this confirms that these two curves are isomorphic to $\PP^1$.
As expected, the row for $p=11$ coincides with the number of points on the elliptic curve 121b1 given above.

\section{Examples}\label{ex_sec}
Throughout this final section, we will give elliptic curves with their label from Cremona's tables~\cite{cremona} together with links to the corresponding page at~\cite{lmfdb}.
The rational \textsc{cm}~points, i.e., points in $X(\QQ)$ whose elliptic curves are defined over $\QQ$ and have complex multiplication, appear frequently; they will be denoted by $E_{D}$ where $D$ is the discriminant of the order.
The list of all of them was given in Table~\ref{cm_overQ_tab}.

\subsection{A necklace on an elliptic curve over \texorpdfstring{$\QQ$}{Q}}
Among the easiest examples to present is the necklace on the \href{https://www.lmfdb.org/EllipticCurve/Q/864a1/}{curve
$E\colon y^2 = x^3 - 3\, x + 6$} for~$p=3$.
The $3$-isogeny field is given by
\[
  K = \QQ[s] / \bigl(s^8 - 6\,s^5 + 3\,s^4 + 18\,s^3 + 18\,s^2 + 18\,s + 9 \bigr)
\]
whose Galois group is $D_4$ over~$\QQ$.
It is \href{https://www.lmfdb.org/NumberField/8.0.2985984.1}{listed in the \textsc{lmfdb}}.
The four subgroups of order $3$ are defined by the following kernel polynomials:
\begin{align*}
  C_0 \colon \ x &-  \tfrac{12}{11}s^7 +\tfrac{ 26}{33}s^6 - \tfrac{14}{33}s^5 + \tfrac{76}{11}s^4 - \tfrac{94}{11}s^3 - 14s^2 - \tfrac{104}{11}s - \tfrac{105}{11}=0,\\
  C_1\colon\ x &+ \tfrac{2}{33}s^7 - \tfrac{2}{33}s^6 - \tfrac{2}{33}s^5 - \tfrac{8}{33}s^4 + \tfrac{10}{11}s^3 + \tfrac{10}{11}s^2 -\tfrac{14}{11}s - \tfrac{13}{11}=0,\\
  C_2\colon\ x&+\tfrac{ 8}{11}s^7 -\tfrac{ 8}{11}s^6 +\tfrac{ 20}{33}s^5 -\tfrac{ 54}{11}s^4 +\tfrac{ 76}{11}s^3 +\tfrac{ 76}{11}s^2 +\tfrac{ 74}{11}s +\tfrac{ 75}{11}=0,\\
  C_3\colon\ x &+\tfrac{ 10}{33}s^7 -\tfrac{ 4}{33}s^5-\tfrac{ 58}{33}s^4 +\tfrac{ 8}{11}s^3 +\tfrac{ 68}{11}s^2 + 4s +\tfrac{ 43}{11}=0
\end{align*}
The necklace $(C_0,C_1,C_2,C_3)$ is the only necklace defined over~$\QQ$.
The corresponding $j$-invariants are
\begin{align*}
  j(E/C_1) &=\tfrac{1}{11} \bigl( 218818080 s^7 - 67867696s^6 - 146242512s^5 - 1250596464s^4  \\
  &\phantom{=\tfrac{1}{11}\cdot \bigl(} + 1173069456s^3 + 4543203600s^2 + 1837676880s - 318777768\bigr),\\
  j(E/C_2) &= \tfrac{1}{11}\bigl(-3174080s^7 - 33133776s^6 + 68459664s^5 + 34515312s^4 \\
  &\phantom{=\tfrac{1}{11}\cdot \bigl(} - 112554000s^3 - 113068176s^2 - 130229424s - 75810888\bigr),\\
  j(E/C_3) &= \tfrac{1}{11} \bigl( 5891808s^7 + 30416048s^6 + 33887568s^5 + 14692080s^4 \\
  &\phantom{=\tfrac{1}{11}\cdot \bigl(} - 161874960s^3 - 161360784s^2 - 152352720s - 97934184\bigr),\\
  j(E/C_4) &= \tfrac{1}{11}\bigl( -221535808s^7 + 70585424s^6 + 43895280s^5 + 1201389072s^4 \\
  &\phantom{=\tfrac{1}{11}\cdot \bigl(} - 898640496s^3 - 4268774640s^2 - 1555094736s - 3711549384\bigr).
\end{align*}
It is difficult to present examples with larger~$p$.
Typically the coefficients of the polynomial defining the $p$-isogeny field and the $j$-invariants (or the kernel polynomials) have very large height.
The code file at~\cite{implementation} contains more complicated examples.

\subsection{Reduction to a supersingular curve without extra automorphisms}\label{ss_13_5_ex_subsec}
This example is for $p=5$ and we will consider the reduction modulo $\ell=13$ of necklaces on curves with complex multiplication.
The \textsc{cm}~points correponding to the curves $E_{-7}$, $E_{-8}$, $E_{-28}$, $E_{-67}$, and $E_{-163}$ as listed above all have reduction at $\ell=13$ isomorphic to the unique supersingular curve
\[
  \tilde E\colon y^2 = x^3 + x + 4
\]
whose $j$-invariant is~$5\in\FF_{13}$.
This curve has no extra automorphisms.
All six pearls are defined over the quadratic extension $\FF_{169} = \FF_{13}[\theta]/(\theta^2-\theta+2)$.
Here we list the polynomials that define these subgroups on the above model of~$\tilde{E}$.
\begin{align*}
  C_0\colon \ & x^2 + (9+5\theta)\,x + 8+12\theta&\qquad
  \bar{C}_0 \colon\  & x^2 + (1+8\theta)\,x + 7+\theta \\
  C_1 \colon\  & x^2 + (8+7\theta)\,x + 11+2\theta&\qquad
  \bar{C}_1\colon\ & x^2 + (2+6\theta)\,x + 11\theta\\
  C_2 \colon \ & x^2 + (7+5\theta)\,x + 1+6\theta&\qquad
  \bar{C}_2\colon\ & x^2 + (12+8\theta)\,x + 7+7\theta
\end{align*}
Here the bar denotes the conjugate over $\FF_{13}$.
The three necklaces on $\tilde{E}$ pictured below are defined over $\FF_{13}$.
\begin{center}
  \includegraphics[height=32mm]{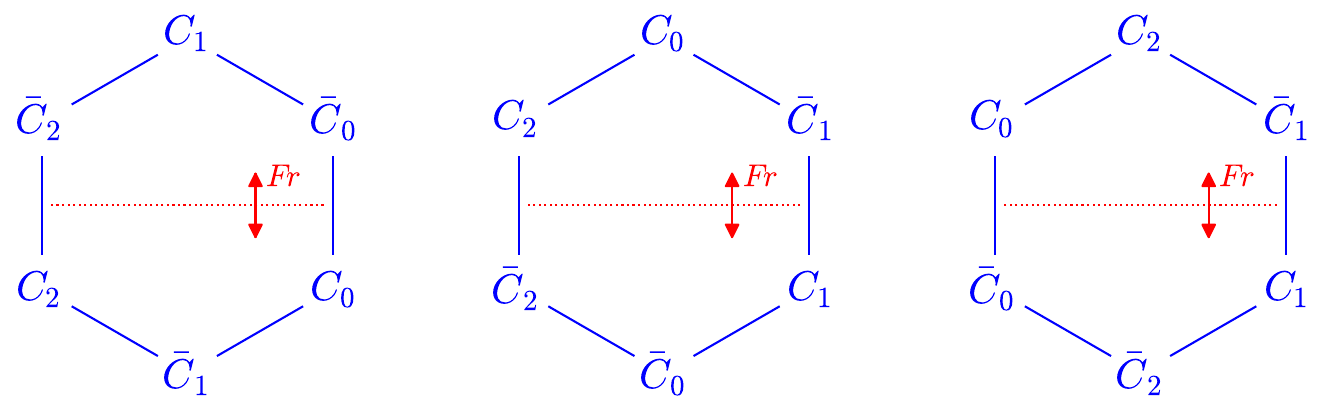}
\end{center}
The pictures are arranged such that the Galois action of $\FF_{169}/\FF_{13}$ is the reflection with respect to the horizontal symmetry.
The unique necklace on the curve $E_{-7}$ reduces to the left-hand necklace.
For the curves $E_{-28}$ and $E_{-67}$, it reduces to the middle necklace, while  the necklaces of $E_{-8}$ and $E_{-163}$ have the right-hand necklace as their reduction.

Note that, further than the three necklaces above, the curve $\tilde E$ has one more $\FF_{13}$-rational necklace, namely $(C_0, C_1,C_2,\bar{C}_0,\bar{C}_1,\bar{C}_2)$.
Frobenius acts on it by rotation.
No curve defined over $\QQ$ with complex multiplication reduces to that necklace.
The remaining six necklaces form three pairs of conjugate necklaces defined over~$\FF_{13^2}$.

\subsection{Reduction to a curve with extra automorphism of order 4}
Now we consider again $p=5$ but $\ell=7$.
The six \textsc{cm}~points with curves $E_{-7}$, $E_{-8}$, $E_{-28}$, $E_{-43}$, $E_{-67}$, and $E_{-163}$
all reduce to the curve $\tilde{E}\colon \ y^2= x^3+x$ with $j=6=1728\in\FF_7$.
This curve has an automorphism $[i]$ of order~4.
The pearls are defined over $\FF_{49} = \FF_7[\theta]/(\theta^2-\theta+3)$.
\begin{align*}
  C_0 & \colon\ x^2 + 6\theta & \qquad
  \bar{C}_0 &\colon\  x^2 + 6 +\theta\\
  C_1 &\colon\   x^2 + 2\theta\, x + 3\theta & \qquad
  \bar{C}_1 &\colon\  x^2 + (2+5\theta)\,x + 3+ 4\theta \\
  C_2 &\colon\  x^2 + 5\theta\, x + 3\theta & \qquad
  \bar{C}_2&\colon\  x^2 + (5+2\theta)\,x + 3+4\theta
\end{align*}
The action by the Frobenius of $\FF_{49}/\FF_7$ is indicated by the bar.
The extra automorphism acts as an involution:
\[
  [i](C_0) = C_0,\ [i](C_1) = C_2,\ [i](\bar{C}_1) = \bar{C}_2.
\]
The left necklace below is the reduction of $E_{-28}$; it is fixed by both $[i]$ and the Frobenius and hence it represents on its own a point of~$X(\FF_7)$.
The picture on the right is the reduced necklace for $E_{-7}$ and $E_{-163}$, which is distinct from $E_{-28}$, but has the same sort of action by Galois and the automorphisms.

\begin{center}
  \includegraphics[height=32mm]{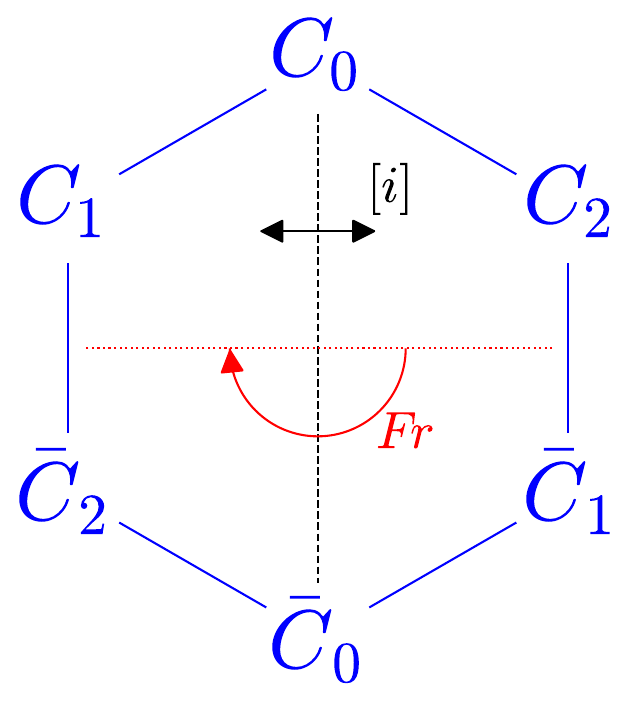}
  \hspace{1cm}
  \includegraphics[height=32mm]{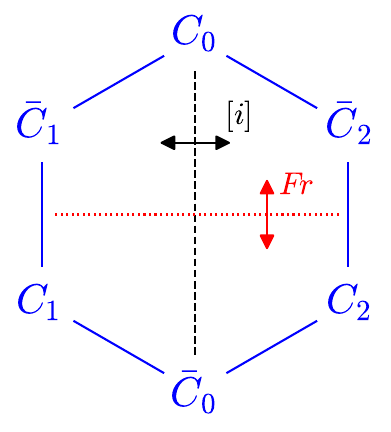}
\end{center}
The next picture is the reduction of~$E_{-8}$.
The point in the modular curve is represented by a pair of necklaces exchanged by $[i]$.
Frobenius exchanges the two necklaces, which shows that the point is $\FF_7$-rational.

\begin{center}
  \includegraphics[height=32mm]{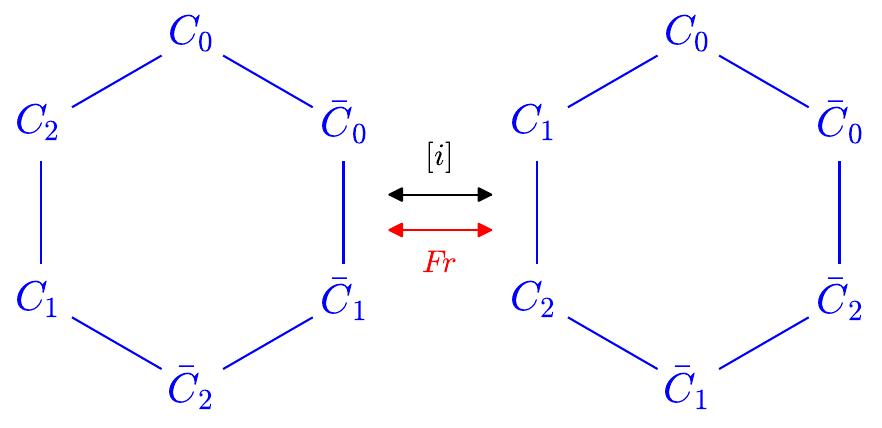}
\end{center}
The final picture is the reduction of both $E_{-43}$ and~$E_{-67}$.
Here the point in $X(\FF_7)$ is again formed by a pair of necklaces exchanged by~$[i]$.
Instead each necklace is already defined over~$\FF_7$.
\begin{center}
  \includegraphics[height=32mm]{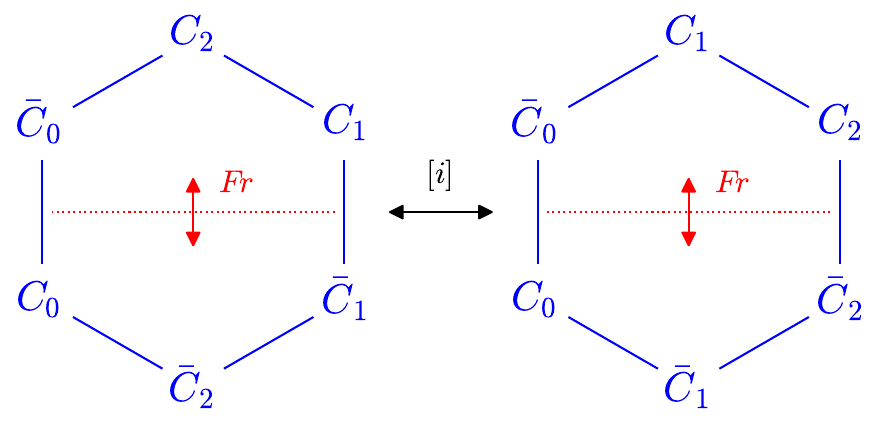}
\end{center}
These four are all the $\FF_7$-rational points on $X$, the remaining four necklaces produce two conjugate points in $X(\FF_{49})$.

\subsection{Reduction to a curve with extra automorphism of order 6}
In this example, we consider again $p=5$, but now $\ell=11$, and we concentrate on the curve $\tilde{E}\colon y^2 = x^3+1$ with~$j=0$.
It has extra endomorphisms and we denote $[\zeta]$ one of the elements of order $3$.
Two pearls are defined over $\FF_{11}$, while the other four are defined over $\FF_{121}= \FF_{11}[\theta]/(\theta^2 + 7\theta + 2)$:
\begin{align*}
  C_0&\colon\   x^2 + 5\,x + 1 &\qquad
  C_1&\colon\  x^2 + 7\,x + 8\\
  C_2&\colon\  x^2 + (4+5\theta)\,x + 7+10\theta   & \qquad
  \bar{C}_2&\colon  x^2 + (2+6\theta)\,x + 3+\theta\\
  C_3&\colon\ x^2 + (10+7\theta)\,x + 1+ 3\theta & \qquad
  \bar{C}_3&\colon  x^2 + (5+4\theta)\,x + 2+8\theta
\end{align*}
The action of the automorphisms satisfies
\[
  [\zeta](C_0) = C_2,\ \ [\zeta](C_2) = \bar{C}_2,\ \ [\zeta](C_1) = C_3, \ \text{ and }\ [\zeta](C_3)=\bar{C}_3.
\]
The reduction of $E_{-3}$ gives the following necklace:
\begin{center}
  \includegraphics[height=32mm]{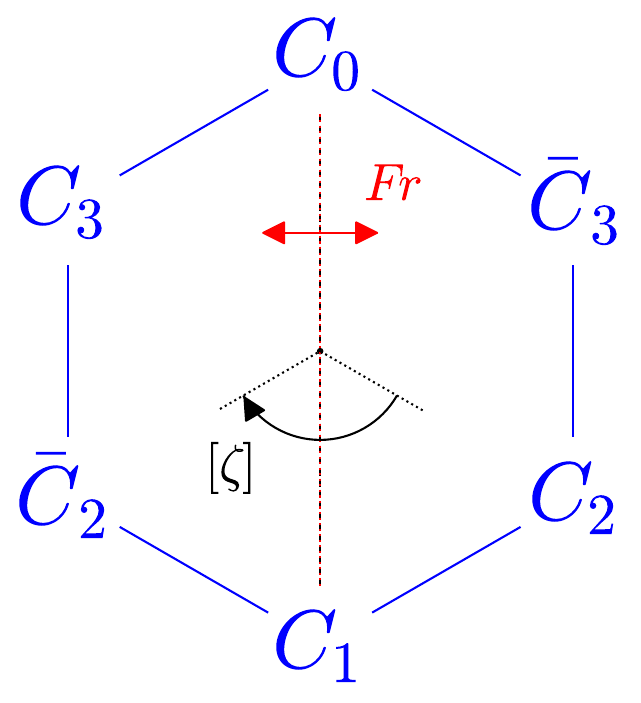}
\end{center}
This necklace is flipped by Frobenius and fixed by $[\zeta]$; therefore it represents a point in~$X(\FF_{11})$.

Instead, the reduction of $E_{-67}$ is the point in $X(\FF_{11})$ represented by the triple of necklaces in the following picture:
\begin{center}
  \includegraphics[height=32mm]{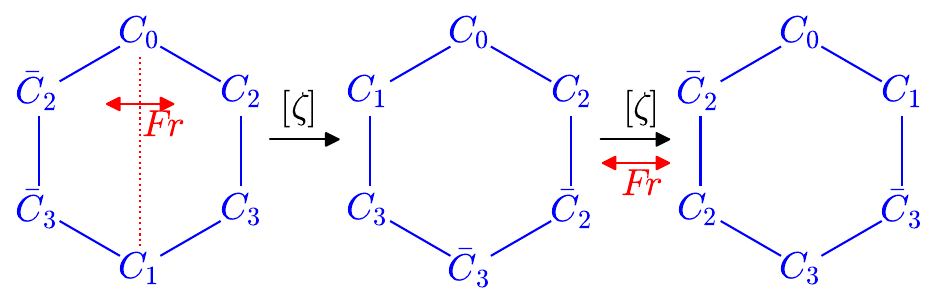}
\end{center}
While the first necklace is fixed by Frobenius, the other two are exchanged by it.

These are the only two $\FF_{11}$-rational points on $X$, the remaining six necklaces on $\tilde E$ form two conjugate points in $X(\FF_{121})$.
While no other rational elliptic curves with complex multiplication that reduce to~$\tilde{E}$, there are of course plenty of curves without complex multiplication.

\subsection{Examples of a necklace with \texorpdfstring{$p=7$}{p=7}}
There are four \textsc{cm}~points corresponding to the curves $E_{-8}$, $E_{-11}$, $E_{-67}$, and $E_{-163}$,
that reduce to the supersingular curve with $j=5$ modulo~$\ell=13$.

As in Section~\ref{ss_13_5_ex_subsec}, the pearls are defined over $\FF_{169} = \FF_{13}[\theta]/(\theta^2-\theta+2)$:
\begin{align*}
  C_0\colon\, &x^3 + 10\,x^2 + 11\,x & \
  C_1\colon\, & x^3 + 7\,x^2 + 5\,x + 12  \\
  C_2\colon\, &x^3 + (4+11\theta)x^2 + (11+9\theta)x + 2+3\theta & \
  \bar{C}_2\colon\, & x^3 + (2+2\theta)x^2 + (7+4\theta)x + 5+10\theta \\
  C_3\colon \, & x^3 + (10+4\theta)x^2 + (8+10\theta)x + 3+12\theta & \
  \bar{C}_3\colon \, & x^3 + (1+9\theta)x^2 + (5+3\theta)x + 2+\theta\\
  C_4 \colon \, &x^3 + 9\,x^2 + (8+\theta)x + 10+8\theta &\
  \bar{C}_4\colon\, & x^3 + 9\,x^2 + (9+12\theta)x + 5+ 5\theta
\end{align*}
There are three necklaces defined over $\FF_{13}$ on this curve:
\begin{center}
  \includegraphics[height=32mm]{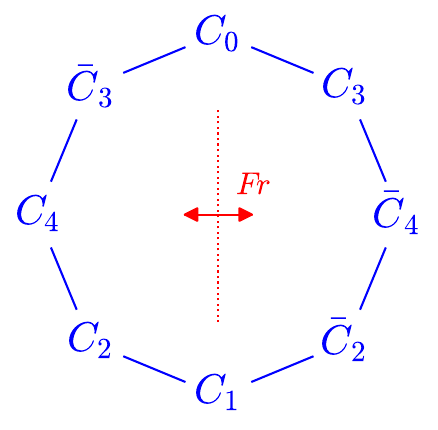}
  \hspace{7mm}
  \includegraphics[height=32mm]{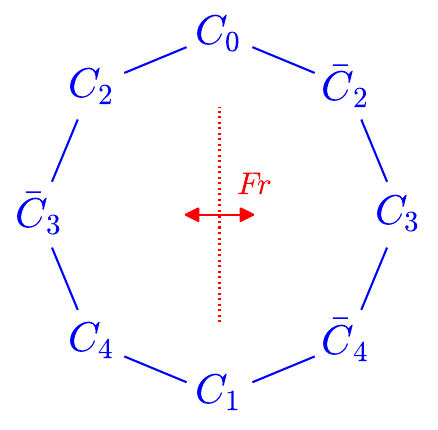}
  \hspace{7mm}
  \includegraphics[height=32mm]{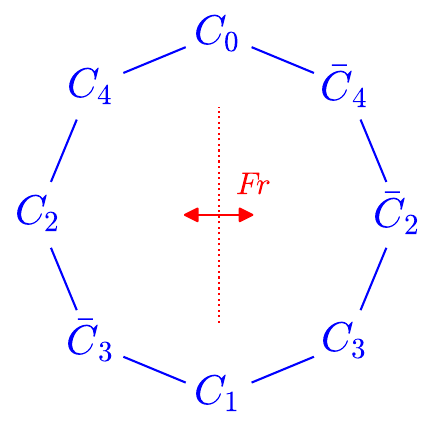}
\end{center}
The left is the reduction of the necklace on $E_{-11}$, the middle is the reduction of the necklace on $E_{-8}$, while the right necklace is the reduction of the necklaces of both $E_{-67}$ and~$E_{-163}$.

\appendix
\section{Appendix: Number of points in the reduction}\label{nu_pts_app}

Here is a table that allows for a simple algorithm to count the number of points in $X(\FF_{\ell})$ for a prime number $\ell\neq p$.

For $j\in\PP^1(\FF_p)$, we write $X_j$ for the fibre of $X(\FF_{\ell})\to \PP^1(\FF_{\ell})$ above~$j$.
Therefore $X(\FF_{\ell}) = \bigcup_{j\in\PP^1(\FF_{\ell})} X_j$.

The following table determines $\# X_j$ in all cases.
Here $a \in \FF_p$ is the reduction modulo $p$ of the trace of Frobenius of the elliptic curve $E$ with the corresponding $j$-invariant.
We define $\delta = \bigl(\tfrac{a^2-4\ell}{p}\bigr)\in\{-1,0,1\}$ and $i$ to denote the number of isogenies on $E$ of degree $p$ defined over $\FF_{\ell}$.
The last invariant can take longer to calculate, but we need it only rarely.

The cases for $j\not\in\{0,1728\}$ at the bottom of the table were obtained in Lemma~\ref{number_of_necklaces_lem}.
the other entries can be calculated by analysing the Galois action on elliptic and ramified points in the special fibres above $j=0$ and $j=1728$.
The calculations are fairly simple, but too tedious to reproduce here.

\begin{center}
\begin{tabular}{clc}
  \toprule
  $j$ & conditions & $\#X_j$ \\
  \midrule
  $\infty$ & $\ell\equiv \pm 1 \pmod{p}$ & $(p-1)/2$ \\
  $\infty$ & $\ell\not \equiv \pm 1 \pmod{p}$ & $0$ \\
  $0$ & $\ell\equiv 2 \pmod{3}$ and $\delta =1$ & $(p-1)/2$ \\
  $0$ & $\ell\equiv 2 \pmod{3}$ and $\delta =-1$ & $(p+3)/2$ \\
  & $\ell\equiv p\equiv 1 \pmod{3}$ and $a=0$ & $(p-1)/6$ \\
  & $\ell\equiv p \equiv 1 \pmod{3}$ and $a^2\equiv 3\ell\pmod{p} $ & $(p-1)/6$ \\
  & $\ell\equiv p \equiv 1 \pmod{3}$ and $a^2\equiv \ell \pmod{p} $ & $p(p-1)/6$ \\
  & $\ell\equiv p \equiv 1 \pmod{3}$ and $a^2\equiv 4\ell\pmod{p} $ & $p(p-1)/6$ \\
  & $\ell\equiv p \equiv 1 \pmod{3}$ and not above & $0$ \\
  & $\ell\equiv 1 \pmod{3}$ and $p \equiv 2\pmod{3}$ and $a=0$ & $(p+7)/6$ \\
  & $\ell\equiv 1 \pmod{3}$ and $p \equiv 2\pmod{3}$ and $a^2\equiv 3\ell\pmod{p}$ & $(p+7)/6$ \\
  & $\ell\equiv 1 \pmod{3}$ and $p \equiv 2\pmod{3}$ and $a^2\equiv \ell\pmod{p}$ & $(p^2-p+4)/6$ \\
  & $\ell\equiv 1 \pmod{3}$ and $p \equiv 2\pmod{3}$ and $a^2\equiv 4\ell\pmod{p}$ & $(p^2-p+4)/6$ \\
   $0$  & $\ell\equiv 1 \pmod{3}$ and $p \equiv 2\pmod{3}$ and not above & $1$ \\
  $1728$ & $p\equiv 1 \pmod{4}$ and $a\neq 0$ and $\delta=0$ & $(p^2-1)/4$\\
  $1728$ & $p\equiv 1 \pmod{4}$ and $a\neq 0$ and $\delta=1$ & $0$ \\
   & $p\equiv 1 \pmod{4}$ and $a= 0$ and $\delta=-1$ & $(p+3)/2$ \\
   & $p\equiv \ell\equiv 1 \pmod{4}$ and $a= 0$ and $\delta=1$ & $(p^2-1)/4$ \\
   & $p\equiv 1 \pmod{4}$ and not above & $(p-1)/2$ \\
   & $p\equiv 3 \pmod{4}$ and $a\neq 0$ and $\delta=0$ & $(p^2+3)/4$ \\
   & $p\equiv 3 \pmod{4}$ and $a\neq 0$ and $\delta\neq 0$ & $1$ \\
  & $p\equiv 3 \pmod{4}$ and $a= 0$ and $\delta = 1$ & $(p-1)/2$ \\
  & $p\equiv\ell\equiv 3 \pmod{4}$ and $a= 0$ and $\delta - 1$ & $(p+3)/4$ \\
  $1728$ & $p\equiv 3 \pmod{4}$ and not above & $(p^2+3)/4$\\
  others & $a=0$ and $\delta =1 $ & $(p-1)/2$ \\
  others & $a=0$ and $\delta =-1$ & $(p+3)/2$ \\
  & $a\neq 0$ and $\delta = 1$ & $0$ \\
  & $a\neq 0$ and $\delta=-1$ & $1$ \\
  & $a\neq 0$, $\delta=0$ and $i>2$ & $p(p-1)/2$ \\
  others & $a\neq 0$, $\delta=0$ and $i\leq 1$ & $0$ \\
  \bottomrule
\end{tabular}
\end{center}

\DeclareEmphSequence{\itshape}
\bibliographystyle{amsplain}
\bibliography{cnl}

@article {cursed,
    AUTHOR = {Balakrishnan, Jennifer and Dogra, Netan and M\"uller, J.
              Steffen and Tuitman, Jan and Vonk, Jan},
     TITLE = {Explicit {C}habauty-{K}im for the split {C}artan modular curve
              of level 13},
   JOURNAL = {Ann. of Math. (2)},
  FJOURNAL = {Annals of Mathematics. Second Series},
    VOLUME = {189},
      YEAR = {2019},
    NUMBER = {3},
     PAGES = {885--944},
}

@article {cursed2,
    AUTHOR = {Balakrishnan, Jennifer S. and Dogra, Netan and M\"uller, J.
              Steffen and Tuitman, Jan and Vonk, Jan},
     TITLE = {Quadratic {C}habauty for modular curves: algorithms and
              examples},
   JOURNAL = {Compos. Math.},
  FJOURNAL = {Compositio Mathematica},
    VOLUME = {159},
      YEAR = {2023},
    NUMBER = {6},
     PAGES = {1111--1152},
       DOI = {10.1112/s0010437x23007170},
       URL = {https://doi.org/10.1112/s0010437x23007170},
}

@article {baran,
    AUTHOR = {Baran, Burcu},
     TITLE = {Normalizers of non-split {C}artan subgroups, modular curves, and the class number one problem},
   JOURNAL = {J. Number Theory},
  FJOURNAL = {Journal of Number Theory},
    VOLUME = {130},
      YEAR = {2010},
    NUMBER = {12},
     PAGES = {2753--2772}
}

@article {geometric_chabauty,
    AUTHOR = {Edixhoven, Bas and Lido, Guido},
     TITLE = {Geometric quadratic {C}habauty},
   JOURNAL = {J. Inst. Math. Jussieu},
  FJOURNAL = {Journal of the Institute of Mathematics of Jussieu. JIMJ.
              Journal de l'Institut de Math\'ematiques de Jussieu},
    VOLUME = {22},
      YEAR = {2023},
    NUMBER = {1},
     PAGES = {279--333},
       DOI = {10.1017/S1474748021000244},
       URL = {https://doi.org/10.1017/S1474748021000244},
}

@article {broker_lauter_sutherland,
    AUTHOR = {Br\"{o}ker, Reinier and Lauter, Kristin and Sutherland, Andrew
              V.},
     TITLE = {Modular polynomials via isogeny volcanoes},
   JOURNAL = {Math. Comp.},
  FJOURNAL = {Mathematics of Computation},
    VOLUME = {81},
      YEAR = {2012},
    NUMBER = {278},
     PAGES = {1201--1231},
       DOI = {10.1090/S0025-5718-2011-02508-1},
       URL = {https://doi.org/10.1090/S0025-5718-2011-02508-1},
      NOTE = {the files we use can be found at \url{https://math.mit.edu/~drew/ClassicalModPolys.html}}
}

@article {campagna-pengo,
  AUTHOR = {Campagna, Francesco and Pengo, Riccardo},
   TITLE = {Entanglement in the family of division fields of elliptic
            curves with complex multiplication},
 JOURNAL = {Pacific J. Math.},
 FJOURNAL = {Pacific Journal of Mathematics},
   VOLUME = {317},
    YEAR = {2022},
  NUMBER = {1},
   PAGES = {21--66},
   ISSN = {0030-8730,1945-5844},
    DOI = {10.2140/pjm.2022.317.21},
     URL = {https://doi.org/10.2140/pjm.2022.317.21},
}

@article {bourdon_clark,
    AUTHOR = {Bourdon, Abbey and Clark, Pete L.},
     TITLE = {Torsion points and {G}alois representations on {CM} elliptic
              curves},
   JOURNAL = {Pacific J. Math.},
  FJOURNAL = {Pacific Journal of Mathematics},
    VOLUME = {305},
      YEAR = {2020},
    NUMBER = {1},
     PAGES = {43--88},
       DOI = {10.2140/pjm.2020.305.43},
       URL = {https://doi.org/10.2140/pjm.2020.305.43},
}

@book {cremona,
    AUTHOR = {Cremona, John E.},
     TITLE = {Algorithms for modular elliptic curves},
   EDITION = {Second},
 PUBLISHER = {Cambridge University Press},
      YEAR = {1997}
}

@misc{furio_lombardo,
      title={Serre's uniformity question and proper subgroups of ${C}_{ns}^+(p)$},
      author={Lorenzo Furio and Davide Lombardo},
      year={2023},
      eprint={2305.17780},
      archivePrefix={arXiv},
      primaryClass={math.NT},
      note={available at \url{https://arxiv.org/abs/2305.17780}},
}

@article {gross_zagier,
    AUTHOR = {Gross, Benedict H. and Zagier, Don B.},
     TITLE = {On singular moduli},
   JOURNAL = {J. Reine Angew. Math.},
  FJOURNAL = {Journal f\"ur die Reine und Angewandte Mathematik. [Crelle's
              Journal]},
    VOLUME = {355},
      YEAR = {1985},
     PAGES = {191--220},
       DOI = {10.1515/crll.1985.355.191},
       URL = {https://doi.org/10.1515/crll.1985.355.191},
}

@article {halberstadt,
    AUTHOR = {Halberstadt, Emmanuel},
     TITLE = {Sur la courbe modulaire {$X_{\text{nd\'ep}}(11)$}},
   JOURNAL = {Experiment. Math.},
  FJOURNAL = {Experimental Mathematics},
    VOLUME = {7},
      YEAR = {1998},
    NUMBER = {2},
     PAGES = {163--174}
}

@book {lang,
    AUTHOR = {Lang, Serge},
     TITLE = {Elliptic functions},
    SERIES = {Graduate Texts in Mathematics},
    VOLUME = {112},
   EDITION = {Second},
      NOTE = {With an appendix by J. Tate},
 PUBLISHER = {Springer-Verlag, New York},
      YEAR = {1987},
     PAGES = {xii+326},
       DOI = {10.1007/978-1-4612-4752-4},
       URL = {https://doi.org/10.1007/978-1-4612-4752-4},
}

@article {lefourn_lemos,
    AUTHOR = {Le Fourn, Samuel and Lemos, Pedro},
     TITLE = {Residual {G}alois representations of elliptic curves with
              image contained in the normaliser of a nonsplit {C}artan},
   JOURNAL = {Algebra Number Theory},
  FJOURNAL = {Algebra \& Number Theory},
    VOLUME = {15},
      YEAR = {2021},
    NUMBER = {3},
     PAGES = {747--771}
}

@incollection {ligozat,
    AUTHOR = {Ligozat, G{\'e}rard},
     TITLE = {Courbes modulaires de niveau {$11$}},
 BOOKTITLE = {Modular functions of one variable, {V} ({P}roc. {S}econd {I}nternat. {C}onf., {U}niv. {B}onn, {B}onn, 1976)},
     PAGES = {149--237. Lecture Notes in Math., Vol. 601},
 PUBLISHER = {Springer},
   ADDRESS = {Berlin},
      YEAR = {1977}
}

@misc{lmfdb,
  shorthand    = {LMFDB},
  author       = {The {LMFDB Collaboration}},
  title        = {The {L}-functions and modular forms database},
  howpublished = {\url{https://www.lmfdb.org}},
  year         = {2025},
  note         = {information about modular curves are currently only available at \url{https://beta.lmfdb.org/ModularCurve/Q/?family=Xnsplus}},
}

@misc{love_studnia_vonk,
      title={Arithmetic intersections on non-split Cartan modular curves},
      author={Jonathan Love and Elie Studnia and Jan Vonk},
      year={2026},
      eprint={2604.06963},
      archivePrefix={arXiv},
      primaryClass={math.NT},
      url={https://arxiv.org/abs/2604.06963},
      note={available at \url{https://arxiv.org/abs/2604.06963}},
}

@article {mercuri_schoof,
    AUTHOR = {Mercuri, Pietro and Schoof, Ren\'e},
     TITLE = {Modular forms invariant under non-split {C}artan subgroups},
   JOURNAL = {Math. Comp.},
  FJOURNAL = {Mathematics of Computation},
    VOLUME = {89},
      YEAR = {2020},
    NUMBER = {324},
     PAGES = {1969--1991},
       DOI = {10.1090/mcom/3503},
       URL = {https://doi.org/10.1090/mcom/3503},
}

@article {necklaces,
    AUTHOR = {Rebolledo, Marusia and Wuthrich, Christian},
     TITLE = {A moduli interpretation for the non-split {C}artan modular
              curve},
   JOURNAL = {Glasg. Math. J.},
  FJOURNAL = {Glasgow Mathematical Journal},
    VOLUME = {60},
      YEAR = {2018},
    NUMBER = {2},
     PAGES = {411--434}
}

@misc{implementation,
  author = {Rebolledo, Marusia and Wuthrich, Christian},
  title = {{SageMath} implementation of necklaces on elliptic curves},
  year = {2025},
  howpublished = {available at \url{https://www.maths.nottingham.ac.uk/plp/pmzcw/}}
}

@incollection {rubin,
    AUTHOR = {Rubin, Karl},
     TITLE = {Elliptic curves with complex multiplication and the conjecture
              of {B}irch and {S}winnerton-{D}yer},
 BOOKTITLE = {Arithmetic theory of elliptic curves ({C}etraro, 1997)},
    SERIES = {Lecture Notes in Math.},
    VOLUME = {1716},
     PAGES = {167--234},
 PUBLISHER = {Springer, Berlin},
      YEAR = {1999}
}

@book {serre,
    AUTHOR = {Serre, Jean-Pierre},
     TITLE = {Lectures on the {M}ordell-{W}eil theorem},
    SERIES = {Aspects of Mathematics},
   EDITION = {Third},
      NOTE = {Translated from the French and edited by Martin Brown from notes by Michel Waldschmidt, With a foreword by Brown and Serre},
 PUBLISHER = {Friedr. Vieweg \& Sohn},
   ADDRESS = {Braunschweig},
      YEAR = {1997},
     PAGES = {x+218}
}

@article {serre72,
    AUTHOR = {Serre, Jean-Pierre},
     TITLE = {Propri\'et\'es galoisiennes des points d'ordre fini des
              courbes elliptiques},
   JOURNAL = {Invent. Math.},
  FJOURNAL = {Inventiones Mathematicae},
    VOLUME = {15},
      YEAR = {1972},
    NUMBER = {4},
     PAGES = {259--331},
      ISSN = {0020-9910,1432-1297},
       DOI = {10.1007/BF01405086},
       URL = {https://doi.org/10.1007/BF01405086}
}

@inproceedings {sutherland_volcanoes,
    AUTHOR = {Sutherland, Andrew V.},
     TITLE = {Isogeny volcanoes},
 BOOKTITLE = {A{NTS} {X}---{P}roceedings of the {T}enth {A}lgorithmic
              {N}umber {T}heory {S}ymposium},
    SERIES = {Open Book Ser.},
    VOLUME = {1},
     PAGES = {507--530},
 PUBLISHER = {Math. Sci. Publ., Berkeley, CA},
      YEAR = {2013}
}

@manual{sagemath,
  Key          = {SageMath},
  Author       = {{The Sage Developers}},
  Title        = {{S}age{M}ath, the {S}age {M}athematics {S}oftware {S}ystem ({V}ersion 10.5)},
  note         = {\url{https://www.sagemath.org}},
  Year         = {2025},
}

@misc{zywina2015,
      title={On the possible images of the mod $\ell$ representations associated to elliptic curves over $\mathbb {Q}$},
      author={David Zywina},
      year={2015},
      eprint={1508.07660},
      archivePrefix={arXiv},
      primaryClass={math.NT},
      url={https://arxiv.org/abs/1508.07660},
      note={available at \url{https://arxiv.org/abs/1508.07660}}
}

\end{document}